\documentclass[10pt]{article}

 \usepackage{grffile}
 \usepackage{ctable}

 \usepackage[small]{caption}
 \usepackage{graphicx}

\usepackage{amsmath,amsfonts,amssymb}
 \usepackage{amsthm}
\usepackage{color}
 \usepackage[letterpaper,textheight=9in,textwidth=7in]{geometry}

 \usepackage[skip=10pt,font=footnotesize]{caption}
\graphicspath{{figures/}}

 \newtheorem{Lemma}{Lemma}[section]
 \newtheorem{Theorem}{Theorem}

  {\begin{trivlist} \item[]{\bf Proof. }}%
  {\hspace*{\fill}$\rule{.4\baselineskip}{.4\baselineskip}$\end{trivlist}}
 
  {\begin{trivlist}\item[]\textbf{Acknowledgments.}}{\end{trivlist}}

  {\begin{center}\textbf{Abstract}}{\end{center}}

 \makeatletter\@addtoreset{figure}{section}\makeatother

 \makeatletter \@addtoreset{equation}{section} \makeatother

\def\Xint#1{\mathchoice
{\XXint\displaystyle\textstyle{#1}}%
{\XXint\textstyle\scriptstyle{#1}}%
{\XXint\scriptstyle\scriptscriptstyle{#1}}%
{\XXint\scriptscriptstyle\scriptscriptstyle{#1}}%
\!\int}
\def\XXint#1#2#3{{\setbox0=\hbox{$#1{#2#3}{\int}$ }
\vcenter{\hbox{$#2#3$ }}\kern-.6\wd0}}

\def\dashint{\Xint-}
\newcommand{\R}{\mathbb{R}}

\newcommand{\Z}{\mathbb{Z}}

\newcommand{\rmi}{\mathrm{i}}
\newcommand{\rmd}{\mathrm{d}}
\newcommand{\rme}{\mathrm{e}}

\newcommand{\rmO}{\mathcal{O}}
\renewcommand{\Re}{\mathrm{Re}\,}

\newcommand{\eps}{{\varepsilon}}
\begin{document}
%\maketitle
 \begin{center}
 {\fontsize{15}{15}\fontfamily{cmr}\fontseries{b}\selectfont{Strain and defects in oblique stripe growth}}\\[0.2in]
 Kelly Chen$\,^1$, Zachary Deiman$\,^2$, Ryan Goh$\,^3$, Sally Jankovic$\,^2$ and Arnd Scheel$\,^2$ \\[0.1in]
 \textit{\footnotesize 
 $\,^1$Massachusetts Institute of Technology, Department of Mathematics, 182 Memorial Drive, Cambridge, MA 02139,  USA\\
 $\,^2$University of Minnesota, School of Mathematics,   206 Church St. S.E., Minneapolis, MN 55455, USA\\[-0.04in]
 $\,^3$Boston University,  Department of Mathematics and Statistics,    111 Cummington Mall, Boston, MA 02215, USA}
 \date{\small \today} 
 \end{center}

\begin{abstract}
We study stripe formation in two-dimensional systems under directional quenching in a phase-diffusion approximation including non-adiabatic boundary effects. We find stripe formation through simple traveling waves for all angles relative to the quenching line using an analytic continuation procedure. We also present comprehensive analytical asymptotic formulas in limiting cases of small and large angles as well as small and large quenching rates. Of particular interest is a regime of small angle and slow quenching rate which is well described by the glide motion of a boundary dislocation along the quenching line. A delocalization bifurcation of this dislocation leads to a sharp decrease of strain created in the growth process at small angles. We complement our results with numerical continuation reliant on a boundary-integral formulation. We also compare results in the phase-diffusion approximation numerically to quenched stripe formation in an anisotropic Swift Hohenberg equation.
\end{abstract}

%\begin{keywords}
%striped phase, Swift-Hohenberg, phase diffusion, dislocation, defect, directional quenching
%\end{keywords}
%\begin{AMS}
% 35B36, 37C29, 35B32
%\end{AMS}
% 
% \setlength{\parskip}{4pt}
% \setlength{\parindent}{0pt}

\section{Introduction}

% talk more about stripes as a crystalline phase, homogenized description as phase diffusion -> equation for the strain as nabla phi, how can boundary conditions be homogenized: effective selected strain through inhomogeneous Neumann, with d_x depending on d_y and speed

% obtain those effective boundary conditions through computation of coherent structures at the microscale, resolving non-adiabatic effects

% obtain some universal laws for those homogenized boundary conditions from asymptotics

% tradeoff: accurate through computation of laws in SH directly, more universality from computation and analysis in PD

% zero speed, oblique angle, boundary conditions inhomogeneous Neumann, select minimum energy strain

% zero speed, zero angle: boundary conditions select strain dependent on microphase

% positive speed, zero angle: selected strain decreases rapidly, beekie paper

% positive speed, small angle: selected strain non-monotone

We investigate the influence of boundary conditions on the formation of striped patterns. Striped patterns occur in many experimental setups \cite{zigzageutectic,eshel,bodenschatz,icicle,double,bradley,digit,sand,thomas,dipstripe} and their existence and stability is quite well studied. In particular, idealized periodic striped patterns in unbounded, planar systems occur in families parameterized by the wavenumber, the orientation, and a phase encoding translations. Stability depends only on the wavenumber and instability mechanisms include Eckhaus and zigzag instabilities. Away from instabilities, striped phases are well described by a phase diffusion equation for a phase $\varphi$ which encodes the (local) shift of a fixed reference pattern. Local wavenumbers and orientation are encoded in the gradient $\nabla \varphi$. Rigorous derivations are possible in a slow modulation approximation \cite{dsss}. In a homogeneously quenched  pattern-forming system, posed with small noisy initial conditions, the observed pattern indeed locally resembles a suitably rotated and stretched periodic pattern, away from isolated points or lines where defects form. More regular patterns emerge when the pattern-forming region expands in time, either through apical growth at the boundary of the domain, or through directional quenching where a parameter in the system is changed spatio-temporally such that the parameter region where pattern formation is enabled grows temporally. Our interest here is with this growth scencario in an idealized situation. 

A prototypical model equation for the the formation of striped patterns is the Swift-Hohenberg equation 
\begin{equation}\label{e:SH}
 u_t=-(\Delta_{x,y}+1)^2u+\mu u - u^3,\quad (x,y)\in\R^2, \ u\in\R, 
\end{equation}
which, for $\mu>0$, possesses families of stable periodic  striped patterns given through  $u_\mathrm{per}(k x;k)=u_\mathrm{per}(kx+2\pi;k)$, close to $\sqrt{4\mu/3}\cos(kx)$ for small $\mu$ and $k\sim 1$. Directional quenching here refers to the situation where $\mu=-\mu_0\, \mathrm{sign}\,(x-c_xt)$ for some $\mu_0\gtrsim 0$. For patterns with trivial $y$-dependence and $c_x=0$, there exists a family of ``quenched'' periodic patterns $u$ with 
\begin{equation}\label{e:1dsd}
 |u(x)-u_\mathrm{per}(k_xx-\varphi;k_x)|\to 0,\ x\to -\infty,\qquad |u(x)|\to 0, \ x\to +\infty,
\end{equation}
for wavenumbers obeying the strain-displacement relation $k_x=g(\varphi)\sim 1+\frac{\mu_0}{16}\sin\varphi$; see \cite{morrissey,scheelweinburd2017}. 

For positive speeds $c_x>0$, one observes the formation of stripes with a selected wavenumber. This stripe formation is enabled by time-periodic solutions $u(t,x)=u_*(x-c_x t,k_x x)$, with $u_*(\xi,\zeta)=u_*(\xi,\zeta+2\pi)$ and 
\[
 u_*(\xi,\zeta)\to u_\mathrm{per}(\zeta;k_x), \ \xi\to -\infty, \qquad 
 u_*(\xi,\zeta)\to 0, \ \xi\to+\infty.
\]
These solutions represent stripes parallel to the quenching interface $x=c_x t$, with trivial $y$-dependence. The wavenumber $k_x$ of stripes selected by this directional quenching process can be computed in terms of the strain-displacement relation and effective diffusivities $d_\mathrm{eff}$ as 
\[
 k_x\sim k_\mathrm{min}+k_1 c_x^{1/2}+\rmO(c_x^{3/4}), \quad k_1=-\zeta(1/2)\sqrt{2 k_\mathrm{min}/d_\mathrm{eff}},
\] 
where $k_\mathrm{min}$ denotes the minimum of the strain-displacement relation; see \cite{beekie}. 

Including possible $y$-dependence, one would be interested in solutions that create periodic patterns at a given angle relative to the quenching interface. This problem was analyzed in \cite{zigzag} when stripes are nearly perpendicular to the quenching interface and in \cite{gs4} when stripes are almost parallel to the boundary for fixed $c_x>0$. Our focus here is on the case of stripes almost parallel to the quenching interface and small speeds. Most of our results are concerned with a phase-diffusion approximation but we demonstrate numerically good agreement with Swift-Hohenberg computations. 

The phase-diffusion approximation for stripes relies on writing solutions $u$ to \eqref{e:SH} in the form $u(t,x)=u_\mathrm{per}(\varphi;k)$, with $|\nabla_{x,y}\varphi|\sim 1$, slowly varying, and 
\[
 \varphi_t=\Delta \varphi,
 \]
after possibly scaling $x$ and $y$ so that effective diffusivities agree.  Of course, this assumes that the patterns considered here are away from possible instabilities, where for instance the Cross-Newell equations would be more appropriate. 
In a context of directional quenching, such an approximation is meaningful only in the pattern forming region $x<c_x t$. The equation therefore needs to be supplemented at the quenching line $x=c_x t,y\in\R$, with an effective boundary condition, which in particular should reflect the strain-displacement relation in the parallel case with $c_x=0$. We then arrive at 
\begin{equation}\label{e:pd}
 \varphi_t=\Delta \varphi + c_x \varphi_x,\ x<0;\qquad \varphi_x=g(\varphi),\ x=0,
\end{equation}
where $g$ reflects the strain-displacement relation, 
\begin{equation}
 g(\varphi)=g(\varphi+2\pi),\qquad g(\varphi)>0,
\end{equation}
for instance $g(\varphi)=1+\kappa \sin(\varphi)$ for some $0\leq \kappa<1$. Clearly, setting $\varphi=\varphi_*(x)$ and $c_x=0$, we find simple affine profiles for any $\varphi_0\in\R$,
\[
 \varphi_*(x)=\varphi_0+g(\varphi_0)x,
\]
corresponding to the solutions in \eqref{e:1dsd} compatible with the strain-displacement relation. Note that \eqref{e:pd} possesses a gauge symmetry that maps solutions $\varphi(t,x)$ to solutions $\varphi(t,x)+2\pi$, reflecting the periodicity of the underlying periodic pattern that is modulated through $\varphi$. It does not possess a continuous symmetry $\varphi(t,x)\mapsto \varphi(t,x)+\bar{\varphi}$, $\bar{\varphi}\in\R$, which would result in $g\equiv const$ and reflect boundary conditions insensitive to the crystalline microstructure. This latter situation arises at leading order when one derives averaged amplitude or phase equations and one can then think of the presence of a nontrivial flux $g$ as a non-adiabatic effect, not visible in averaged approximations. 

The equation \eqref{e:pd} was analyzed in \cite{beekie} for $y$-independent solutions, deriving in particular universal asymptotics for solutions in the cases $c_x\ll 1$ and $c_x\gg 1$. For $c_x\ll1$, excellent agreement with solutions in \eqref{e:SH} and several other prototypical examples of pattern-forming systems was found, including reaction-diffusion, Ginzburg-Landau, and Cahn-Hilliard equations. For bounded initial conditions and $c_x>0$, solutions eventually become time-periodic up to the gauge symmetry, and converge locally uniformly  to linear profiles for large negative $x$, 
\[
 \varphi(t+T,x)=\varphi(t,x)+2\pi, \quad |\varphi(t,x)-(k_x x-\omega t)|\to 0,\ x\to -\infty, \ \omega=c_xk_x, 
\]
for some $T=\frac{2\pi}{\omega}>0$, for given $g>0$. The existence and stability of such solutions with the minimal, 1:1-resonant period $T=\frac{2\pi}{\omega}$ was established generally in \cite{pauthier}. Here, the resonance refers to the frequency of the periodic solution $2\pi/T$ relative to the frequency of patterns generated in the far field $\omega$. In particular, subharmonic solutions $ 2\pi\ell/T=  \omega$, $\ell>1$, are ruled out. 
\begin{figure}
 \includegraphics[width=1\textwidth]{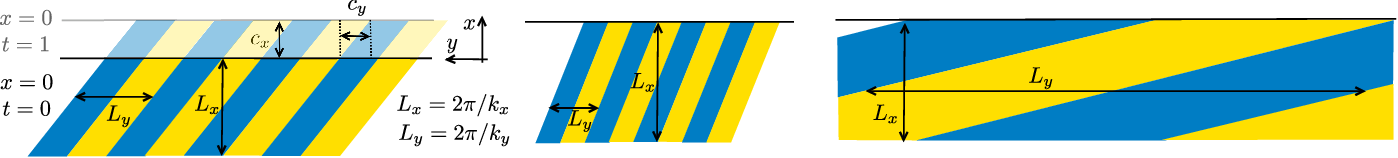}
 \caption{Schematic plot of patterns with found through \eqref{e:pd1}--\eqref{e:pd4}, with $k_y=\rmO(1)$ (left), $k_y\gg1$ (center), and $k_y\ll1$ (right). Also shown on the left is the effect of growth, leading to an apparent drift of the pattern along the interface with speed $c_y=-k_xc_x/k_y$; see text for details. Colors chosen to show contours of $u=u_\mathrm{per}(\varphi(x,y,t))$ with $u_\mathrm{per}(\phi)=\sin(\varphi)$; see Figure \ref{f:profiles} for computed profiles.}\label{f:stsc}
\end{figure}

In the two-dimensional, oblique case, these simplest resonant solutions correspond to traveling waves; see Figure \ref{f:stsc}. In the far field, $x\to -\infty$,  we are interested in oblique stripes which are represented by values of the phase $\varphi\sim k_x (x+c_x t) + k_y y=k_x x+k_y(y-c_y t)$ with $c_y=-k_xc_x/k_y$. Such solutions are in fact traveling waves in the $y$-direction. We therefore  focus on solutions $\varphi(x,k_y(y-c_y t))$ to \eqref{e:pd}, periodic up to the gauge symmetry in the second argument, that is, solutions to 
\begin{alignat}{2}
 0&=\varphi_{xx}+k_y^2\varphi_{\zeta\zeta}+c_x \varphi_x -k_xc_x \varphi_\zeta ,&\quad x<0,\zeta\in\R,\label{e:pd1}\\
 0&=\varphi(x,\zeta+2\pi)-\varphi(x,\zeta)-2\pi,&\quad x\leq 0,\zeta\in\R,\label{e:pd2}\\
 0&=\varphi_x-g(\varphi),&\quad  x=0,\zeta\in\R,\label{e:pd3}\\
 0&=\lim_{x\to-\infty}|\varphi(x,\zeta)-(k_x x + \zeta)|,&\quad \zeta\in\R\label{e:pd4}.
 \end{alignat}
All solutions are in fact classical solutions since we shall assume $g$ to be smooth.  We will also see later that the convergence in \eqref{e:pd4} is in fact uniform. 

In addition to $\varphi$, the system \eqref{e:pd1}--\eqref{e:pd4} includes 3 variables: the lateral periodicity $k_y$, which we will assume to be positive, without loss of generality; the quenching speed $c_x$ which we assume to be non-negative; and the strain $k_x$ in a direction perpendicular to the quenching line, which we think of as a Lagrange multiplier that compensates for the phase shift induced by $\zeta$-translations. Given $k_x=k_x(c_x,k_y)$, one can then determine angle and wavenumber from the wave vector $(k_x,k_y)$. 

Our main results are as follows.

\paragraph{Existence for all $c_x\geq0, k_y>0$}
Assuming $g$ is smooth and $2\pi$-periodic, we have existence.
\begin{Theorem}[Existence]\label{t:ex}
Suppose $g>0$. Then for all $c_x\geq0$, $k_y>0$, we have existence of solutions to \eqref{e:pd1}--\eqref{e:pd4} with $k_x=K_x(k_y,c_x)$, smooth. Moreover, solutions are strictly monotonically increasing in $\zeta$. 
\end{Theorem}
Using reflection symmetry, one can also find monotonically decreasing solutions. Solutions are unique within this class of solutions up to the trivial translation symmetry in $\zeta$.  

\begin{figure}
\centering
\raisebox{0.25in}{\includegraphics[width=0.47\textwidth]{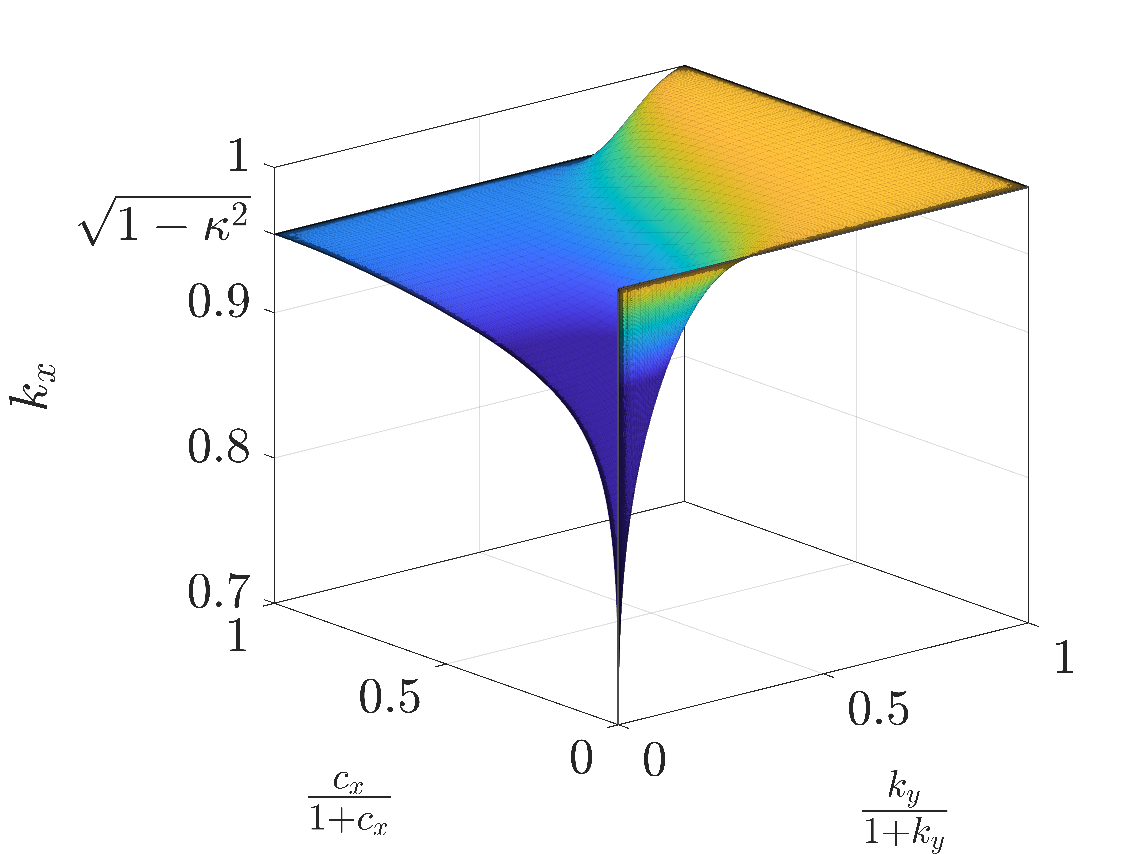}}\hfill
\includegraphics[width=0.51\textwidth]{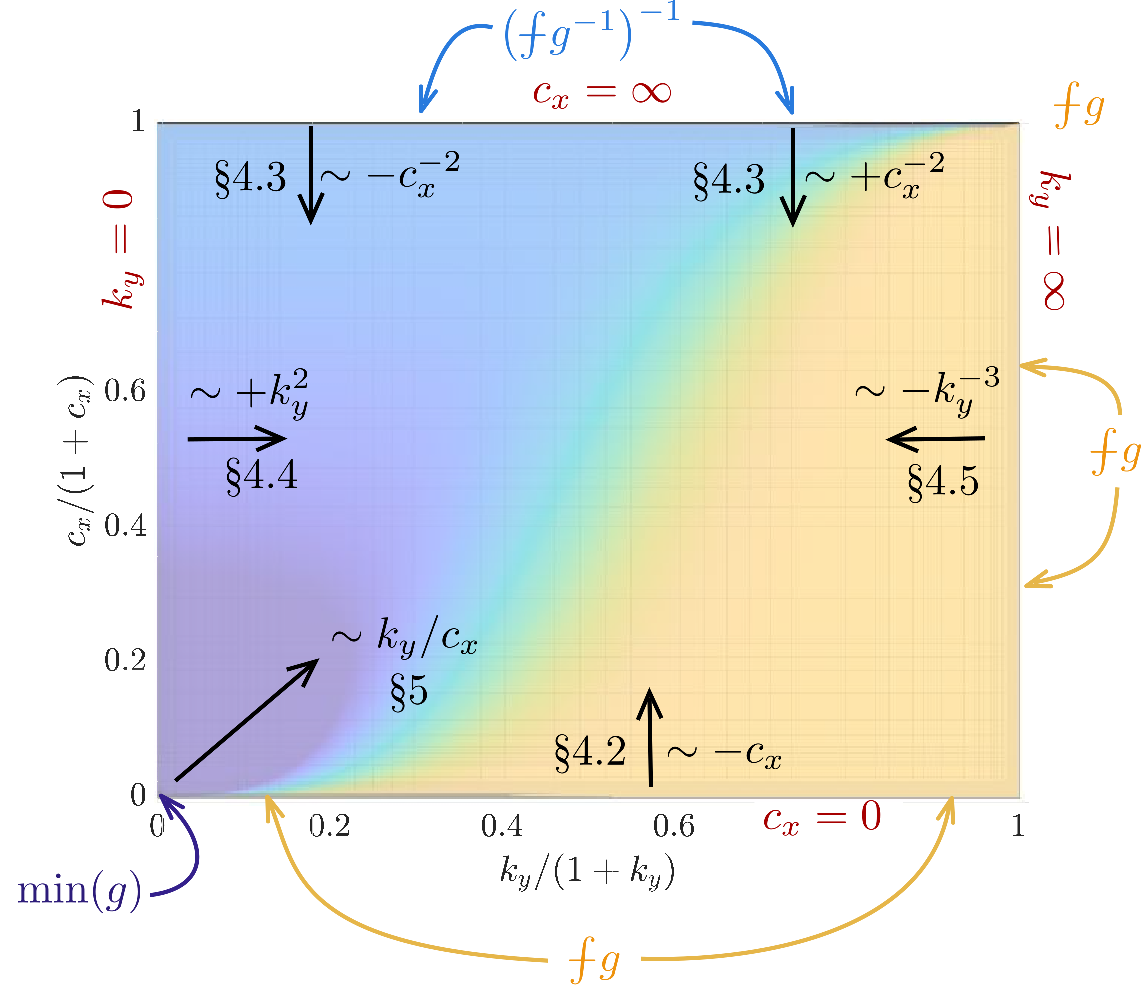}
\caption{Computed values of $k_x$ as a function of $k_y$ and $c_x$ in a compactified scale including the limits $k_y=\infty$ and $c_x=\infty$. Surface plot (left; see \S \ref{s:sum} for other views) and contour plot with limiting values and asymptotics, details in the sections referenced (right).   }\label{f:schemmod}
\end{figure}

We computed the function $K_x(c_x,k_y)$ numerically and show the resulting graph in Figure \ref{f:schemmod}, using an appropriate compactification of the positive quadrant $c_x,k_y\geq 0$. One sees quite distinct limiting behaviors of the surface and much of this paper is concerned with exploring these limits.  Figure \ref{f:schemmod} includes a guide to the asymptotics and how they are reflected in this surface.

\paragraph{Asymptotics $c_x\to\infty$} Solutions $\varphi$ and wavenumbers converge as $c_x\to\infty$ with limiting wavenumber $K_x(c_x=\infty,k_y)$ independent of $k_y$, given through the harmonic average of $g$. At finite but large $c_x$, wavenumbers decrease from the harmonic average for small $k_y$ and increase for large $k_y$, proportional to $c_x^{-2}$ at leading order.

\paragraph{Asymptotics $c_x\to 0$, $k_y>0$ fixed} Solutions and wavenumbers are smooth at $c_x=0$ with limit $k_x$ given by the average of $g$, and linear asymptotics for $c_x$ small. We establish asymptotics for the linear coefficient as $k_y\to 0$.

\paragraph{Asymptotics $k_y\to 0$, $c_x>0$ fixed. } Solutions are smooth (albeit likely not analytic) near $k_y=0$, $c_x>0$, a regime explored also in \cite{gs4}. We numerically compute a leading-order quadratic coefficient and explore asymptotics of this coefficient as $c_x\to 0$ numerically. 

\paragraph{Asymptotics $k_y\to\infty$} In this limit of perpendicular stripes, we find again the average of $g$ as the limit and asymptotics with leading-order term $k_y^{-3}$.

\begin{figure}
\includegraphics[width=0.3\textwidth]{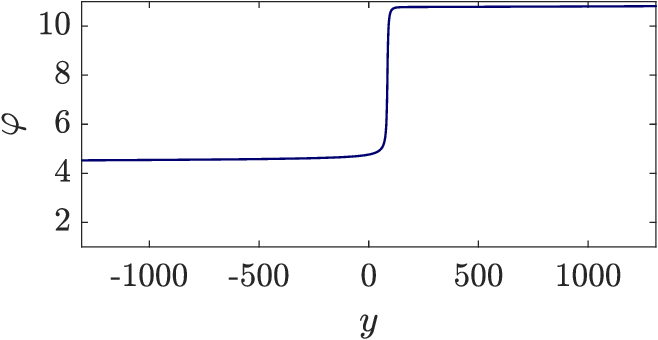}\hfill 
\includegraphics[width=0.68\textwidth]{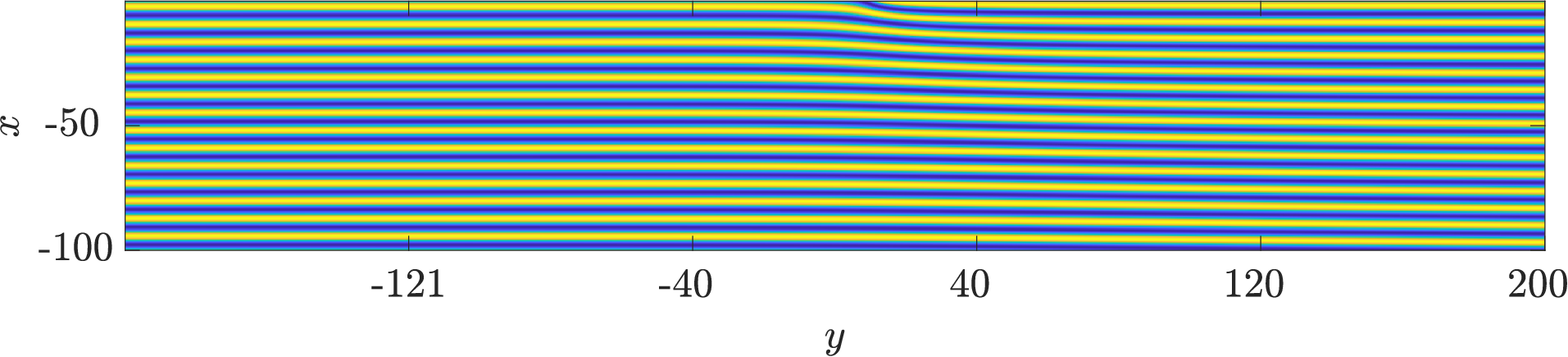}\\[0.2in]
\includegraphics[width=0.3\textwidth]{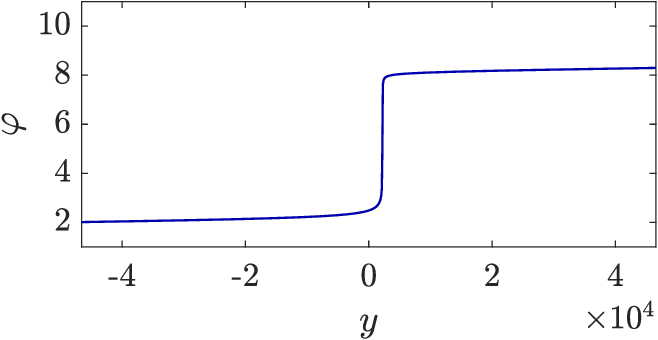}\hfill 
\includegraphics[width=0.68\textwidth]{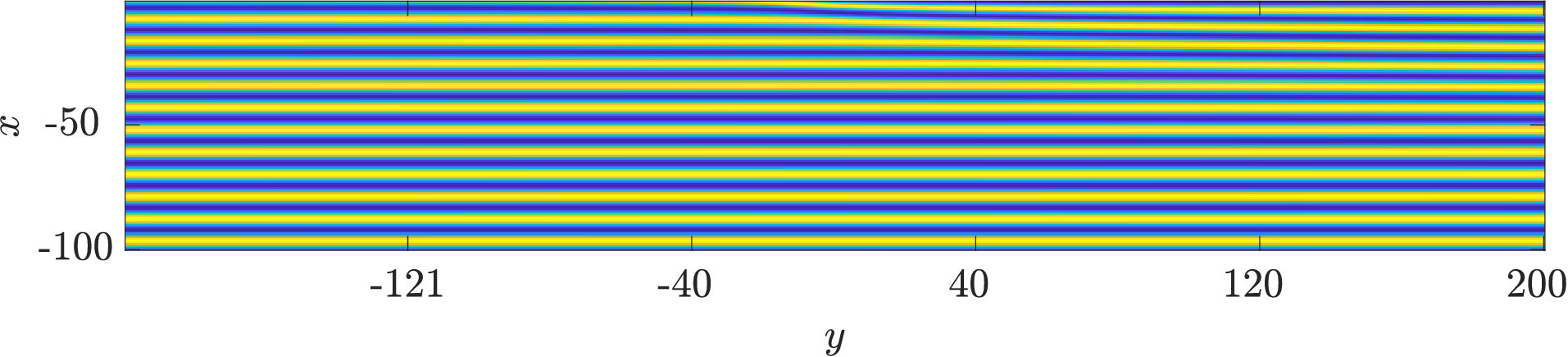}
 \caption{Profiles of $\varphi$ on $x=0$ for $k_y=2.4\times 10^{-3}$, $k_x=0.9$ (top left), and $k_y=6.75\times 10^{-5}$, $k_x=0.7147$ (bottom,left), both with $c_x=10^{-4}$. Note the different scales on the horizontal axis, showing that the jump is stronger localized for larger $k_y$. Associated profiles of $\sin(\varphi)$ in the $x-y$-plane (only part of $y$-region shown), showing a sharply localized defect for larger $k_y$ (top right) and a delocalized defect for small $k_y$ (bottom right).}
 \label{f:profiles}
\end{figure}

\paragraph{Asymptotics $k_y\sim c_x\to 0$}
In the most striking regime close to the origin, the sharp peak in the surface in Figure \ref{f:schemmod}, we use an inner expansion to arrive at a reduced problem which amounts to describing the glide motion of a dislocation-type defect in the $y$-direction under an externally imposed strain. Most interestingly, we identify a qualitative ``phase transition'' where this defect changes type, explaining qualitatively the shape of the surface $k_x(c_x,k_y)$ close to the origin. Profiles of solutions in this regime on the boundary and in the whole domain are shown in Figure \ref{f:profiles}, demonstrating in particular the phase transition corresponding to the delocalization of a defect near $k_y/c_x\sim 2.8845$; see \S\ref{s:origin}.

\paragraph{Numerical continuation} We illustrate results and explore the approximation quality of theoretical asymptotics using numerical continuation for solutions of \eqref{e:pd1}--\eqref{e:pd4}, and also for corresponding solutions of the Swift-Hohenberg equation. We find good agreement with asymptotics in the phase-diffusion equation, and a qualitatively similar transition near $c_x,k_y\sim 0$ due to defect delocalization  in the Swift-Hohenberg equation.

\paragraph{Consequences for homogenized descriptions} Thinking of the gradient of the phase as a macroscopic, homogenized strain variable for a crystalline phase, our results provide corresponding \emph{effective} boundary conditions through a micropscopic analysis of the boundary layer. The dependence $k_x= K_x(k_y;c_x)$ provides mixed boundary conditions, such that the renormalized strain $\phi=\varphi-K_x(k_y;c_x)x$ solves
\[
 \phi_t=\Delta \phi+c_x\phi_x,\ x<0,\qquad \phi_x=0,\ x=0,
\]
eliminating variations on the microscopic scale $1/K_x$. Such a description is not possible for $c_x=k_y=0$, since the derivative $\varphi_x$ at the boundary depends on the microscopic phase variable $\varphi$  and, at steady-state, there are multiple compatible equilibrium strain configurations. The presence of a spatial defect, $k_y\neq 0$, or a temporal defect, $c_x\neq 0$, forces selection of a unique normal strain at the boundary and allows this macroscopic description. From this perspective, our work establishes existence of a unique normal strain and analyzes in detail properties of this normal strain in various limiting regimes, in particular relying on properties of the spatio-temporal defect at the boundary.  We emphasize that these effective boundary conditions are \emph{not} the natural boundary conditions associated with minimizing a free energy density and ``select'' non-energy-minimizing strains; see \S\ref{s:7} for a discussion of stored energies in the growth processes considered here and more context for wavenumber selection in striped phases.

\paragraph{Outline}
We introduce a boundary integral formulation together with a priori estimates and numerical setup in \S\ref{s:2} and prove existence of oblique quenched fronts for all $k_y\neq0, c_x\geq 0$, in \S\ref{s:3}. We derive asymptotics in the limits $c_x\to 0$, $c_x\to\infty$, $k_y\to 0$, and $k_y\to\infty$  in \S\ref{s:4}. We present an analysis near the origin $k_y,c_x\sim 0$ in \S\ref{s:origin} and compare with Swift-Hohenberg in \S\ref{s:anSH}.

% 
% \paragraph{Acknowledgment} KC, ZD, SJ, and AS gratefully acknowledge partial support through grant NSF DMS-1907391.  RG was partially supported through NSF-DMS 2006887.

\section{Boundary integral formulation, a priori estimates, and numerical setup}\label{s:2}
To solve \eqref{e:pd1},\eqref{e:pd2}, and \eqref{e:pd4}, we first set 
\begin{equation}\label{e:aff}
 \psi(x,\zeta):=\varphi(x,\zeta)-(k_x x +\zeta),
\end{equation}
which gives
\begin{alignat}{2}
 0&=\psi_{xx}+k_y^2\psi_{\zeta\zeta}+c_x \psi_x -k_xc_x \psi_\zeta ,&\quad x<0,\zeta\in\R,\label{e:pdp1}\\
 0&=\psi(x,\zeta+2\pi)-\psi(x,\zeta),&\quad x\leq 0,\zeta\in\R, \label{e:pdp2}\\
 0&=\psi_x-g(\psi+\zeta)+k_x,&\quad x=0,\zeta\in\R,\label{e:pdp3}\\
 0&=\lim_{x\to-\infty} \psi(x,\zeta),&\quad \zeta\in\R\label{e:pdp4}.
 \end{alignat}
Next, writing Fourier series $\psi(x,\zeta)=\sum_{\ell\in\Z}\psi^\ell(x) \rme^{\rmi \ell \zeta}$ transforms \eqref{e:pdp1} into
\begin{equation}\label{e:ft}
 \frac{\rmd^2}{\rmd x^2}\psi^\ell+c_x \frac{\rmd}{\rmd x}\psi^\ell -k_y^2\ell^2\psi^\ell -k_x c_x\rmi\ell \psi^\ell=0,
\end{equation}
with 
\begin{equation}
 \psi^\ell(x)=\sum_\pm \psi^\ell_\pm\rme^{\nu_\pm^\ell x},\qquad \nu_\pm^\ell=-\frac{c_x}{2}\pm \sqrt{\frac{c_x^2}{4}+k_y^2\ell^2 +c_x k_x\rmi\ell}, \qquad \ell\neq 0,
\end{equation}
where we use the standard cut at $\R^-$ in the square root and restrict to $c_x\geq 0$. For $\ell\neq 0$, decay \eqref{e:pdp4} requires $\psi_-^\ell=0$. For $c_x=\ell=0$, solutions are affine, $\psi^0(x)=\psi^0_0+\psi^0_1 x$, and we can set $\psi_0^1=0$ since this part of the solution is already parameterized by the ansatz \eqref{e:aff} through the parameter $k_x$.   For $c_x>0$, $\ell=0$, convergence as in  \eqref{e:pdp4} implies $\psi^0(x)\equiv\psi_0^0=0$.
Evaluating $\psi_x$ at $x=0$ and substituting into \eqref{e:pdp3} then reduces \eqref{e:pdp1}--\eqref{e:pdp4} to the boundary-integral equation
\begin{equation}\label{e:bi}
 0= \mathcal{D}_+(\partial_\zeta;c_x,k_x,k_y)\psi - g(\psi+\zeta)+k_x,\quad \psi(\zeta)=\psi(\zeta+2\pi),\quad \mathcal{D}_+(\rmi\ell;c_x,k_x,k_y)=\nu_+^\ell,
\end{equation}
where the operator $\mathcal{D}_{+}$ is understood as a Fourier multiplier acting through multiplication by $\nu_+^\ell$ on Fourier series. One readily confirms that $\mathcal{D}_+:H^1_\mathrm{per}\subset L^2\to L^2$ is a closed, sectorial operator as a relatively compact perturbation of $k_y|\partial_\zeta|$, with compact resolvent and spectrum with strictly positive real part except for the simple eigenvalue $\lambda=0$ associated with constant functions. The definition of $\mathcal{D}_+$ extends to $c_x=0$ in natural agreement with our problem. For later puposes, we also introduce the associated pseudo-differential operator $\mathcal{D}_-$ through $ \mathcal{D}_-(\rmi\ell;c_x,k_x,k_y)=\nu_-^\ell$.

\begin{Lemma}\label{l:apriori}
 For any periodic and smooth flux $g$, there exists  an a priori bound $C_\infty(g,c_x,k_y,m)$ such that any solution to \eqref{e:bi} with $\psi(0)\in [0,2\pi)$ satisfies 
 \[
  \|\psi\|_{C^m}+|k_x|\leq  C_\infty.
 \]
 Moreover, $C_\infty$ is uniformly bounded for fixed $m$ and $\delta>0$ such that $|k_y|>\delta$, $\|g\|_{C^m}\leq 1/\delta$. 
\end{Lemma}
\begin{proof}
 Since the average $\dashint\mathcal{D}_+\psi = 0$, $\dashint =\frac{1}{2\pi}\int$, vanishes and $|g|_\infty\leq C_g$, we find an a priori bound  $|k_x|\leq \dashint |g(\psi(\zeta)+\zeta)|$. This in turn gives an $L^\infty$ a priori bound on $\mathcal{D}_+\psi$ and, using the regularizing properties of $\mathcal{D}_+$ and a bootstrap, the desired a priori bound on $\psi$.  Uniformity of $C_\infty$ follows readily from the fact that the pseudo-inverse of $\mathcal{D}_+$ is uniformly bounded from $L^2$ into $H^{1/2}$ as long as $k_y$ is outside a neighborhood of the origin.
\end{proof}

\paragraph{Numerical setup} We solve \eqref{e:bi} numerically for the variables $\psi$ and $k_x$, with parameters $c_x$ and $k_y$, and adding a phase condition $\int \psi(\zeta)\exp(-\zeta^2/\delta)\rmd \zeta=0$. The resulting nonlinear equation is evaluated using fast Fourier transform. A Newton method, using \texttt{gmres} to solve the linear equation in each Newton step was found to converge robustly even for poor initial guesses. Most of the solutions were then computed using secant continuation in $k_y$ for fixed $c_x$ with adaptive control of the continuation step. During each step, we control for the number of Fourier modes by ensuring that amplitudes in high Fourier modes is below a tolerance, which we found to have little effect once below $10^{-4}$. Step sizes are very small and numbers of Fourier modes grow when $c_x,k_y\sim 0$, due to large gradients in the profile. We address this regime directly using an inner expansion and a slightly different ansatz function in \S\ref{s:origin}. The code was implemented in \texttt{matlab} and Newton iterations for large sizes $N\geq 2^{18}$ were carried out on a Nvidia GV100 GPU. All numerical results use $g(\varphi)=1+\kappa\sin(\varphi)$ with $\kappa=0.3$ unless otherwise noted.

\section{Existence in the phase-diffusion approximation}\label{s:3}

In this section, we prove Theorem~\ref{t:ex}.  Throughout we write $\dashint=\frac{1}{2\pi}\int$ for the average integral. For this, we perform a homotopy, introducing $g_\tau(u):=\tau g(u)+(1-\tau)\dashint g$. Clearly, $g_\tau$ satisfies all the assumptions of Theorem \ref{t:ex} for $\tau\in[0,1]$, in particular $g_\tau>0$. Let $I\subset [0,1]$ be the set of values where the conclusion of Theorem \ref{t:ex} holds. We will show below that 
\begin{center}
(i)\ $0\in I$;\qquad \qquad\qquad
 (ii) $I$ is closed; \qquad \qquad\qquad
 (iii) $I$ is open.
\end{center}
Together, this implies that $I=[0,1]$ and establishes Theorem \ref{t:ex}. This general strategy of proof was used in \cite{pauthier} for the case $k_y=0$, although the proof there was based directly on the parabolic equation rather than the boundary-integral formulation which we shall exploit here. 

To show (i), we set $k_x=\dashint g$ and $\psi=0$, such that $\varphi$ is strictly monotone. 

To show (ii), take a sequence of solutions $\psi^n$ with wavenumbers $k_x^n$ for converging values $\tau^n\to \tau^\infty$. We may assume, possibly adding multiples of $2\pi$, that $\psi^n(0)\in [0,2\pi)$.  By Lemma \ref{l:apriori}, we can assume that $\psi^n\to \psi^\infty$ and $k_x^n\to k_x^\infty$, possibly passing to a subsequence. The limit then solves \eqref{e:pdp1}--\eqref{e:pdp4}. It remains to show that the limit $\varphi^\infty=\psi^\infty+\zeta$ is strictly monotone. Clearly, $(\psi^\infty)'\geq -1$ by uniformity of the limit. We argue by contradiction. Suppose therefore that $(\psi^\infty)'(\zeta_0)=-1$. Note that $v=(\psi^\infty)'+1$ solves \eqref{e:pdp1}, \eqref{e:pdp2}, and \eqref{e:pdp4}, together with the linearized boundary conditions
\[
 0=v_x - g_{\tau^\infty}'(\psi^\infty+\zeta)v, \qquad x=0,\zeta\in\R,
\]
and has $v(\zeta_0)=0$, $v_\zeta(\zeta_0)=0$, $v_{\zeta\zeta}(\zeta_0)\geq 0$. Extending into $x<0$ and using the boundary condition gives $v_x(\zeta_0)=0$ and, using the equation, $v_{xx}\leq 0$. On  the other hand, since $\dashint v >0$ at $x=0$, $v(\zeta,x)\to \dashint v|_{x=0}>0$, a constant. Since interior minima are excluded by the maximum principle, the minimum of $v$ is necessarily located  at the boundary $x=0,\zeta=\zeta_0$, which however implies $v_x(\zeta_0)>0$ by the Hopf boundary lemma, a contradiction. 

It remains to show (iii) for any $g_\tau$. Therefore, first notice that the linearization of \eqref{e:bi} at any profile $\psi_*$,
\[
 \mathcal{L}_* v=\mathcal{D}_+(k_x)v-g_\tau'(\psi_*+\zeta)v,
\]
is Fredholm of index zero with $\psi_*'(\zeta)+1$ belonging to the kernel. We claim that the kernel is indeed one-dimensional and that the derivative of \eqref{e:bi} with respect to $k_x$, $\mathcal{D}'_+(k_x) \psi_*$, does not belong to the range. Together, this then establishes (iii) via the Implicit Function Theorem since the linearization with respect to $(\psi,k)$ is onto. Suppose first that there is a function $v$ in the kernel that is not a multiple of $\psi_*'(\zeta)+1$. Then we can find a linear combination that is non-negative but not strictly positive, that is, a function $w$ in the kernel with $w(\zeta_0)=0$, $w(\zeta)\geq0$, and $\dashint w\geq 0$. Arguing as in (ii), we can then obtain a contradiction from the maximum principle. It now only remains to show that there does not exist a nontrivial solution to 
\begin{equation}\label{e:solv}
 \mathcal{D}_+(k_x) v-g_\tau'(\psi_*+\zeta)v = -\mathcal{D}'_+(k_x) \psi_* -1, 
\end{equation}
where we suppressed the dependence of $\mathcal{D}_+$ on its arguments other than $k_x$. Note that in the case $c_x=0$, $\mathcal{D}_+'=0$, $\mathcal{D}_+$ is self-adjoint, with cokernel $\psi_*'+1$, such that the right-hand side of \eqref{e:solv} has nonzero scalar product with the cokernel and hence does not belong to the range. We shall therefore assume in the sequel that $c_x>0$. The boundary integral equation \eqref{e:solv} is equivalent to the elliptic equation 
\begin{alignat}{2}
 0&=v_{xx}+k_y^2 v_{\zeta\zeta}-c_xk_x v_\zeta + c_x v_x,&\quad x<0,\\
 0&=v_x-g_\tau'(\psi+\zeta)v+\mathcal{D}_+'(k_x)v+1,&\quad x=0.\label{e:gen1}
\end{alignat}
We claim that the existence of a solution to \eqref{e:gen1} is equivalent to the existence of a generalized eigenvector in an associated elliptic problem, which will then lead to a contradiction. Consider therefore the eigenvalue problem associated with our linearization
\begin{alignat}{2}
 0&=v_{xx}+k_y^2 v_{\zeta\zeta}-c_xk_x v_\zeta + c_x v_x -\lambda v,&\quad x<0,\\
 0&=v_x-g_\tau'(\psi+\zeta)v,&\quad x=0,\label{e:gen2a}
\end{alignat}
with solution $v=\psi_*'+1$ at $\lambda=0$. Existence of a generalized eigenvector then amounts to a solution $v$ to 
\begin{alignat}{2}
 0&=v_{xx}+k_y^2 v_{\zeta\zeta}-c_xk_x v_\zeta + c_x v_x - c_x (\psi_*'+1) ,&\quad x<0,\\
 0&=v_x-g_\tau'(\psi+\zeta)v,&\quad x=0,\label{e:gen2b}
\end{alignat}
or, setting $v=w+x$, 
\begin{alignat}{2}
 0&=w_{xx}+k_y^2 w_{\zeta\zeta}-c_xk_x w_\zeta + c_x w_x - c_x \psi_*' ,&\quad x<0,\\
 0&=w_x-g_\tau'(\psi+\zeta)w+1.&\quad x=0.\label{e:gen2c}
\end{alignat}
Solving the first equation using Fourier series in $\zeta$ and a variation-of-constant formula exploiting boundedness as $x\to\infty$,  we find after a short calculation
\[
 w_x(0)=\mathcal{D}_+w(0)+ (\mathcal{D}_+-\mathcal{D}_-)^{-1}c_x \psi_*'|_{x=0},
\]
which is equivalent to \eqref{e:gen1}. This however contradicts the simplicity of the first eigenvalue of the elliptic operator defined in \eqref{e:gen1}.

\section{Asymptotics near the boundaries of $\{k_y>0,c_x>0\}$}\label{s:4}

We derive asymptotics in the regular and singular limits when either $c_x$ or $k_y$ tend to 0 or infinity. 

\subsection{The case $c_x=0$}
In this case, we can multiply \eqref{e:bi} by $\psi_*'+1$ and integrate over $\zeta\in[0,2\pi]$ to find 
\begin{align*}
 0&=\int_\zeta \left((\psi_*'+1)\mathcal{D}_+\psi_* -(\psi_*'+1)g(\psi_*+\zeta) + (\psi_*'+1)k_x\right)\\
  &=2\pi (k_x-\dashint_\varphi g(\varphi) ),
  \end{align*}
where we used that $\mathcal{D}_+$ is a symmetric operator with kernel spanned by the constant functions to see that the first summand vanished, and monotonicity of $\psi_*+\zeta$ to transform the second summand into an integral over $\varphi$. As a consequence $k_x=\dashint g$ is a priori known; see also \cite{lloydscheel,zigzag,bakker}, where this wavenumber selection mechanism was derived from Hamiltonian identities. 

\subsection{The limit $c_x\to 0$}
We suppose that $k_y>0$ and study the limit $c_x\to 0$. Since the operator $\mathcal{D}_+(c_x)$ is continuous in the limit $c_x=0$ as a map from $H^1$ into $L^2$, this limit is a regular perturbation problem. Using in addition that the \begin{figure}[h]
 \centering
 \includegraphics[width=0.45\textwidth]{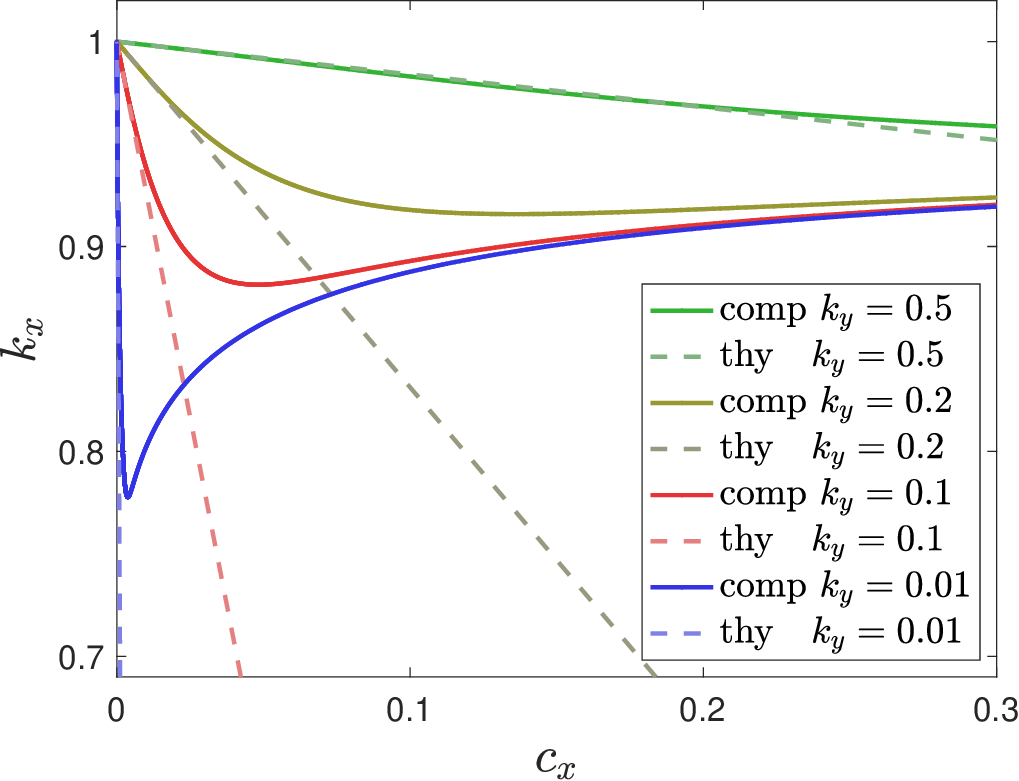}\hfill
 \includegraphics[width=0.445\textwidth]{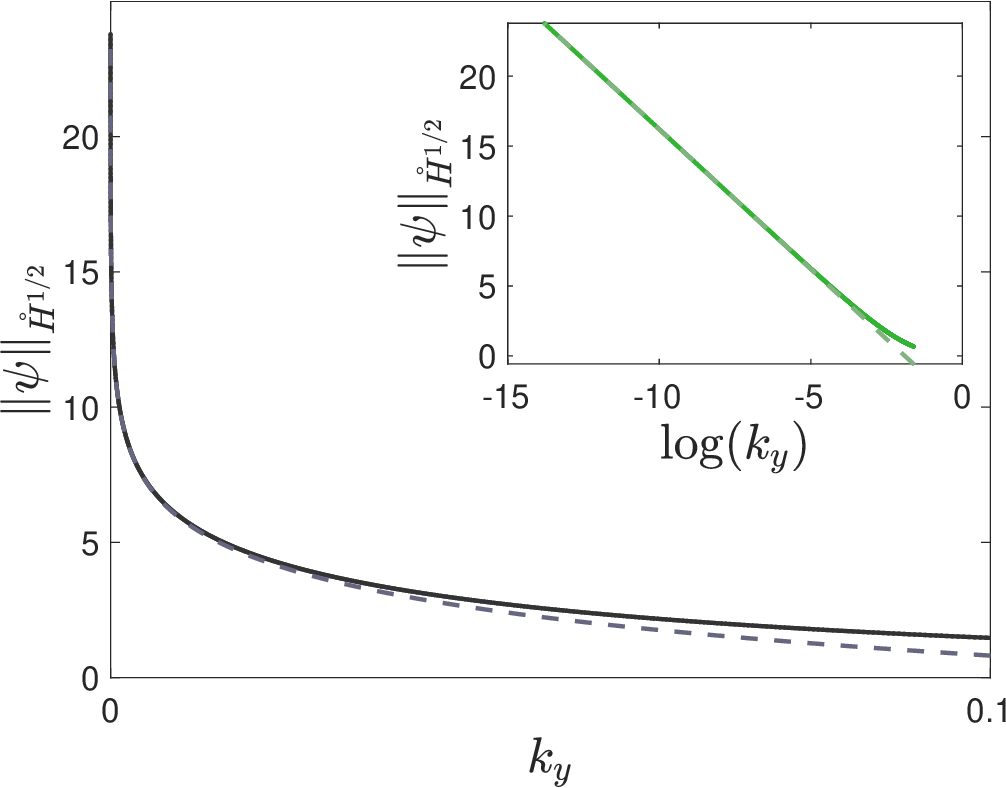}
 \caption{Left: Strain as a function of small growth rates comparing numerical continuation with theory \eqref{e:kx1}, where the $\mathring{H}^{1/2}$-norm was computed numerically. Right: Asymptotics for the $\mathring{H}^{1/2}$-norm as $k_y\to 0$, comparing with theory \eqref{e:h12}, best fit for $\rmO(1)$ terms. }\label{f:1}
\end{figure}
linearization at a profile, including the parameter $k_x$ as a variable, is onto, we conclude that we can formally expand the solution in $c_x$,
\[
 \psi_*(\zeta;c_x)=\psi_*(\zeta;0)+c_x\psi_1(\zeta) + \rmO(c_x^2),\qquad k_x=k_{x,0}+c_x k_{x,1}+\rmO(c_x^2),\ k_{x,0}=\dashint g.
\]
Inserting this expansion into the equation and taking the scalar product with the kernel of the linearization $\psi_*'+1$ gives at order $c_x$, expanding $\mathcal{D}_+=\mathcal{D}_+^0+c_x\mathcal{D}_+^1+\rmO(c_x^2)$, $\mathcal{D}_+^1(\ell)=\frac{1}{2}(-1+\rmi\frac{k_{x,0}}{k_y} \mathrm{sign}(\ell))$, $\ell\neq 0$, $\mathrm{sign}(0)=0$,
\begin{align*}
 0&=\int_\zeta (\psi_*'+1)
 \left(\mathcal{D}_+^1\psi_* + k_{x,1}\right)\\
  &=2\pi \left(k_{x,1}+\frac{k_{x,0}}{2k_y} \|\psi_*\|^2_{\mathring{H}^{1/2}}\right),
\end{align*}
where we used $\int \psi_*=0$ and set $\|\psi_*\|^2_{\mathring{H}^{1/2}}=\dashint \psi|\partial_\zeta|\psi$, which gives
\begin{equation}\label{e:kx1}
 k_x=k_{x,0}+\left(- \frac{k_{x,0}}{k_y}\|\psi_*\|^2_{\mathring{H}^{1/2}}\right)\frac{c_x}{2}+\rmO(c_x^2),\qquad k_{x,0}=\dashint g.
\end{equation}
Formally setting $k_y=c_x=0$, we find that $\mathcal{D}_+=0$ and a solution $\psi_{0,*}=-\zeta\mod 2\pi$ which does not belong to ${\mathring{H}^{1/2}}$. Writing the equation as $\psi = (1+\mathcal{D}_+)^{-1}\psi +g-k_x$ leads to the prediction $\psi_*\sim (1+k_y |\partial_\zeta|)^{-1} \psi_{0,*}$ with 
\begin{equation}\label{e:h12}
 \|\psi_*\|^2_{\mathring{H}^{1/2}} \sim -2 \log(k_y)+\rmO_{k_y}(1).
\end{equation}
In particular, we expect a strong initial stretching , that is, a decrease in $k_x$ with  $c_x$ proportional to $-2c_x|\log(k_y)|/k_y$. 

Computed solutions $k_x$ are compared with the asymptotic prediction in Figure \ref{f:1}, where we also show agreement between the asymptotic prediction for the linear coefficient and the asymptotic formula \eqref{e:h12}.

\subsection{The limit $c_x\to\infty$}
We suppose that $k_y>0$ and study the limit $c_x\to\infty$. We therefore set $c_x=\eps^{-1}$ and formally expand 
\[
\mathcal{D}_+(\ell;\eps)=\rmi k_x \ell + (k_x^2+k_y^2)\ell^2 \eps + (2 \rmi \ell kx (\ell^2 k_x^2 + \ell^2 k_y^2))\eps^2+\rmO(\eps^3).
\]
We start by considering the case $\eps=0$, where $\mathcal{D}_+(\partial_\zeta;0)=k_x\partial_\zeta$. As a consequence, at $\eps=0$, the solution $\psi=\psi^0+\zeta$ solves the ordinary differential equation
\begin{equation}\label{e:0cxinf}
k_x \psi_{0,\zeta}=g(\psi_0),\qquad \psi_0(\zeta+2\pi)=\psi_0(\zeta)+2\pi,
\end{equation}
with implicit solution from separation of variables. In particular, the wavenumber at infinity is the harmonic average of the nonlinearity, 
\[
 k_{x,0}=\left(\dashint \left(g(v)\right)^{-1}\right)^{-1}.
\]
The linearization at $\eps=0$, $\psi^0$ is 
\[
 \mathcal{L}^0 v=k_{x,0} v_\zeta-g'(\psi_0)v,
\]
which we consider as a Fredholm operator of index zero from $H^1_\mathrm{per}$ into $L^2$. The derivative of  \eqref{e:0cxinf} with respect to $k_x$ is $\psi_{0,\zeta}$ which does not belong to the range, so that the linearization is, as in the case of finite $c_x$ discussed in \S\ref{s:3}, onto and we can use the Implicit Function Theorem to solve. Since the equation is not smooth in $\eps$, one needs to be somewhat careful. We therefore first expand formally, 
\[
 k_x=k_{x,0}+k_{x,1}\eps + k_{x,2}\eps^2+\rmO(3),\qquad \psi=\psi_0+\psi_1\eps+\psi_2\eps^2+\rmO(3), 
\]
where $\psi_j,$ $j>1$ are periodic, and substitute into \eqref{e:bi}. At first order, we find 
\begin{equation}\label{e:o1}
\mathcal{L}^0 \psi_1+
 \left( 
   k_{x,1}\psi_{0,\zeta} -
   \left( 
     (k_{x,0})^2 +k_y^2
   \right)
   \psi_{0,\zeta\zeta}  
 \right)=0 .
 \end{equation}
Integrating against the adjoint kernel $1/\psi_{0,\zeta}$ we see that $k_{x,1}=0$ since, using the chain rule to compute $\psi_{0,\zeta\zeta}$ and changing integration to $\psi$ instead of $\zeta$, 
\[
 \int_0^{2\pi} \frac{\psi_{0,\zeta\zeta}}{\psi_{0,\zeta}}\rmd\zeta = \int_0^{2\pi}\frac{g'(\psi)}{g(\psi)}\rmd\psi=0,
 \]
by periodicity of $\log(g(\psi))$. We can then solve for $\psi_1$ as 
\begin{equation}\label{e:phi1}
 \psi_1=\frac{(k_{x,0})^2 +k_y^2}{k_{x,0} }\log(\psi_{0,\zeta})\psi_{0,\zeta}=\frac{(k_{x,0})^2 +k_y^2}{(k_{x,0})^2 }\log\left(\frac{g(\psi_0)}{k_{x,0}}\right)g(\psi_0).
\end{equation}
At order $\eps^2$, we find 
\begin{equation}\label{e:o2}
 \begin{split} \mathcal{L}^0 \psi_2+
  \Big( 
   k_{x,2}\psi_{0,\zeta} -
 \frac{1}{2}g''(\psi_0)(\psi_1)^2 
 +2k_{x,0}  \left( 
     (k_{x,0})^2 +k_y^2
   \right) \psi_{0,\zeta\zeta\zeta} +
   k_{x,1} \psi_{1,\zeta} \\
   -\  2 k_{x,0}k_{x,1} \psi_{0,\zeta\zeta}
   -\left( 
     (k_{x,0})^2 +k_y^2
   \right)\psi_{1,\zeta\zeta}
\Big).
\end{split}
\end{equation}
Using that $k_{x,1}=0$, integrating against the kernel of the adjoint $1/\psi_{0,\zeta}$, and changing variables of integration gives 
\begin{equation*}
k_x=k_{x,0}+ k_{x,2}c_x^{-2}+\rmO(c_x^{-4}),
\end{equation*}
with 
\begin{equation}
 \begin{split}
k_{x,2}=\dashint 
 \Big\{
 \frac{1}{2}g''(\psi_0)(\psi_1)^2 
 -2k_{x,0}  \left( 
     (k_{x,0})^2 +k_y^2
   \right) \psi_{0,\zeta\zeta\zeta} 
   +\left( 
     (k_{x,0})^2 +k_y^2
   \right)\psi_{1,\zeta\zeta}
 \Big\}\\
 \times\frac{1}{(\psi_{0,\zeta})^2}\rmd\psi_0,\label{e:kx2}
\end{split}
\end{equation}
where one substitutes 
\begin{equation}
 \begin{split}
\psi_{0,\zeta}=&\frac{1}{k_{x,0}}g(\psi_0),\\
\psi_{0,\zeta\zeta}=&\frac{1}{(k_{x,0})^2}g'(\psi_0)g(\psi_0),\\
\psi_{0,\zeta\zeta\zeta}=&\frac{1}{(k_{x,0})^3}\left(g''(\psi_0)(g(\psi_0))^2+(g'(\psi_0))^2 g(\psi_0)\right),
\end{split}
\end{equation}
and uses equation \eqref{e:phi1}. 

The resulting integrals can be evaluated numerically for specific choices of $g(v)$. We found that for $g(v)=1+\kappa \sin(v)$, $|\kappa|<1$, $k_{x,2}$ is monotonically increasing as a a function of $k_y$, $k_{x,2}<0$ for $k_y=0$ and $0<k_{x,2}\sim k_y^4$ for $k_y$ large. 
More explicitly, the  integrals can be evaluated to order $\kappa^4$ for $g(\varphi)=1+\kappa\sin(\varphi)$, yielding 
\begin{equation}\label{e:cxinf2kappa}
k_x(c_x)= \sqrt{1-\kappa^2}+ \frac{1}{2}\left(\kappa^2 (-1 + k_y^4) + \frac{1}{4} \kappa^4 (3 + 5 k_y^4)+\rmO(\kappa^6)\right)c_x^2 + \rmO(c_x^4).
\end{equation}
This proves in particular that, at least for small $\kappa$, the monotonicity of $k_x$ as a function of $c_x$ changes, that is, $k_{x,2}$ changes sign, to leading order at $k_y=1$. 

\begin{figure}[h]
 \centering
 \includegraphics[width=0.49\textwidth]{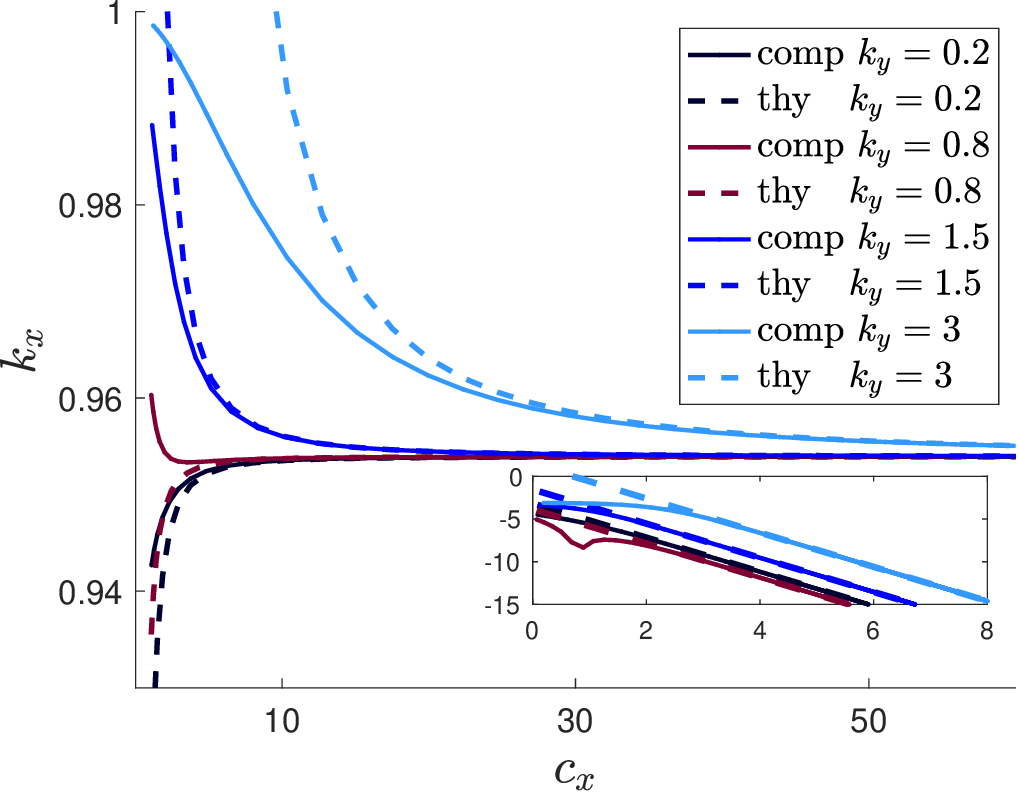}\hfill
 \includegraphics[width=0.475\textwidth]{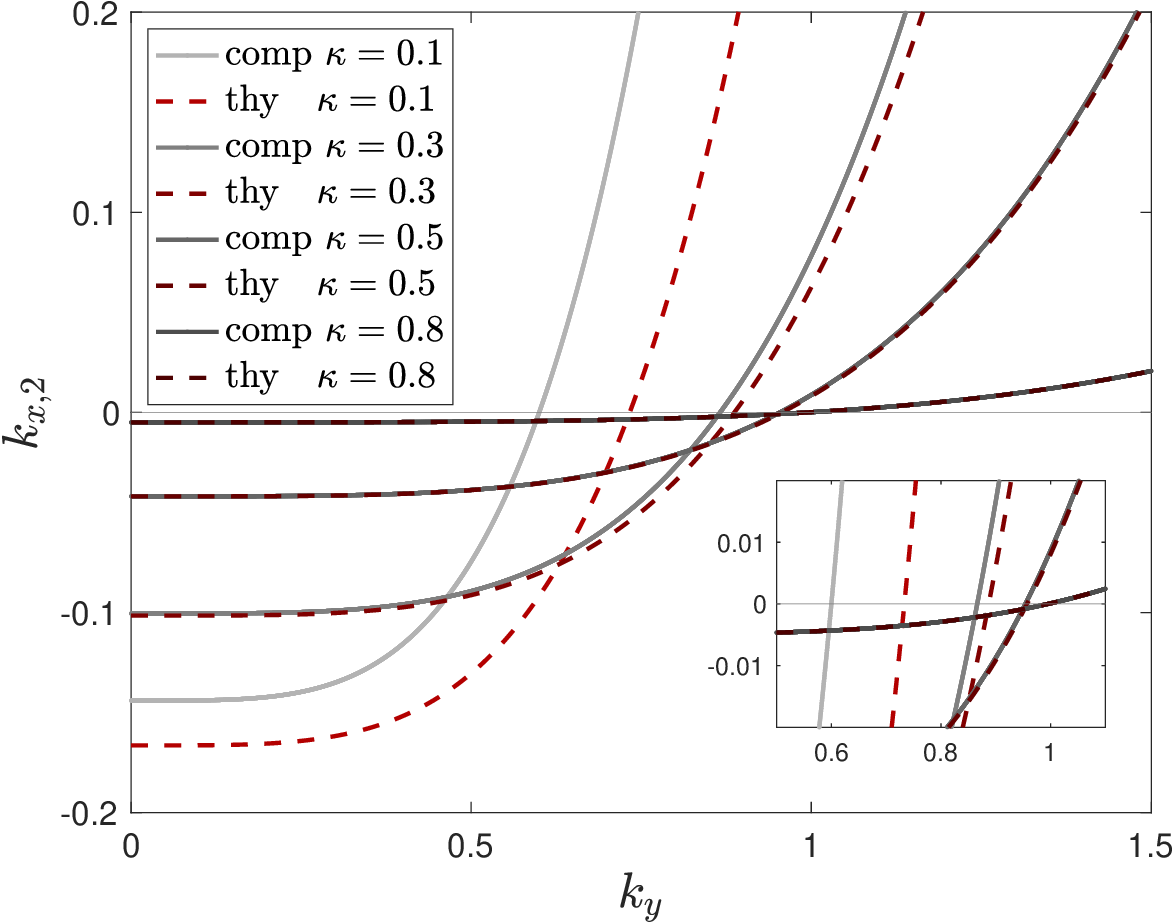}
 \caption{Left: $k_x$ for large $c_x$ for several $k_y$-values, compared with theory \eqref{e:kx2}; inset shows comparison on of $k_x-k_{x,\infty}$ and $c_x$ on $\log$-scales. Right: Leading-order coefficient $k_{x,2}$ as a function of $k_y$ through numerical evaluation of \eqref{e:kx2}, (solid), and explicit approximation \eqref{e:cxinf2kappa},}\label{f:2}
\end{figure}

Figure \ref{f:2} shows numerically computed values of $k_x$ compared with asymptotics for large $c_x$, for several values of $k_y$, and demonstrates the sign change of the second-order coefficient $k_{x,2}$ in a comparison with \eqref{e:cxinf2kappa}.

In order to make this expansion rigorous, we rewrite the equation as 
\begin{equation}\label{e:fac}
 (1-\mathcal{D}_{+,1}(\eps,\zeta))\mathcal{D}_+^0\psi - g(\psi) + k_x=0. 
\end{equation}
The operator $(1-\mathcal{D}_{+,1}(\eps,\zeta))$ is bounded invertible on $L^2$ as a direct inspection of the Fourier symbol shows. Moreover, it is continuous at $\eps=0$ as an operator from $H^1$ to $L^2$, again via a direct inspection of the Fourier symbol, with limit the identity. Therefore, \eqref{e:fac} can be written as 
\[
 F(\psi,k_x):=\mathcal{D}_+^0\psi - (1-\mathcal{D}_{+,1}(\eps,\zeta))^{-1}(g(\psi - k_x)=0,
\]
where $F:H^1\times\R\to L^2$ is continuous in $\eps$ at $\eps=0$. The Implicit Function Theorem then guarantees the existence of solutions for $\eps>0$, small, with leading-order terms $\psi^0,k_{x,0}$. Substituting subsequently higher-order expansion, one can proceed in a similar fashion to establish validity of the expansion to any fixed order. 

\subsection{The limit $k_y\to 0$}
We follow the strategy from the previous section and find 
 at $\rmO(2)$, 
 \[
 0=\dashint \psi^\mathrm{ad}
 \left(
    c_x^2+4 c_x k_{x,0}\partial_\zeta
 \right)^{-1/2}
 \left( 
 \psi_{0,\zeta\zeta}-c_x k_{x,2} \psi_{0,\zeta}
 \right)\rmd\zeta,
 \]
where $\psi^\mathrm{ad}$ is the (unique up to scalar multiples) periodic solution to the adjoint equation 
$\mathcal{D}_+(-\partial_\zeta)\psi_0-g'(\psi_0)\psi_0=0$. Unfortunately, the solution to the adjoint equation does not appear to be readily expressible in terms of $\psi_0$ so that we will rely on numerical methods to evaluate the integral and obtain coefficients $k_{x,0}$ and $k_{x,2}$ in the expansion 
\begin{equation}\label{e:kyto0}
 k_x=k_{x,0}+k_{x,2} k_y^2 + \rmO(k_y^4).
\end{equation}

The numerically computed results shown in Figure \ref{f:3}  show good agreement up to a sharp transition value that we shall discuss in \S\ref{s:origin}. 
 
Numerically, we find that the quadratic coefficient $k_{x,2}$ decreases with $c_x$ in a monotone fashion, converges to $0$ as $c_x\to\infty$ and to $\infty$ for $c_x\to 0$, with power law asymptotics $k_{x,2}\sim c_x^{-\beta}$, $\beta\sim 1/2$. Asymptotics are well captured through
\begin{equation}\label{e:cx2asy}
 k_{x,2}=c_x^{-1/2}(c_1\log(c_x)+c_2);
\end{equation}
fitting $c_1$ and $c_2$ for $c_x\in [5\cdot 10^{-6},1\cdot 10^{-5}]$ provides excellent agreement for a wide range of $c_x$-values; see Figure \ref{f:3}.
We did not attempt to justify asymptotics but provide a conceptual explanation in \S\ref{s:origin}.

\subsection{The limit $k_y\to\infty$}
\begin{figure}[h]
 \centering
 \includegraphics[width=0.39\textwidth]{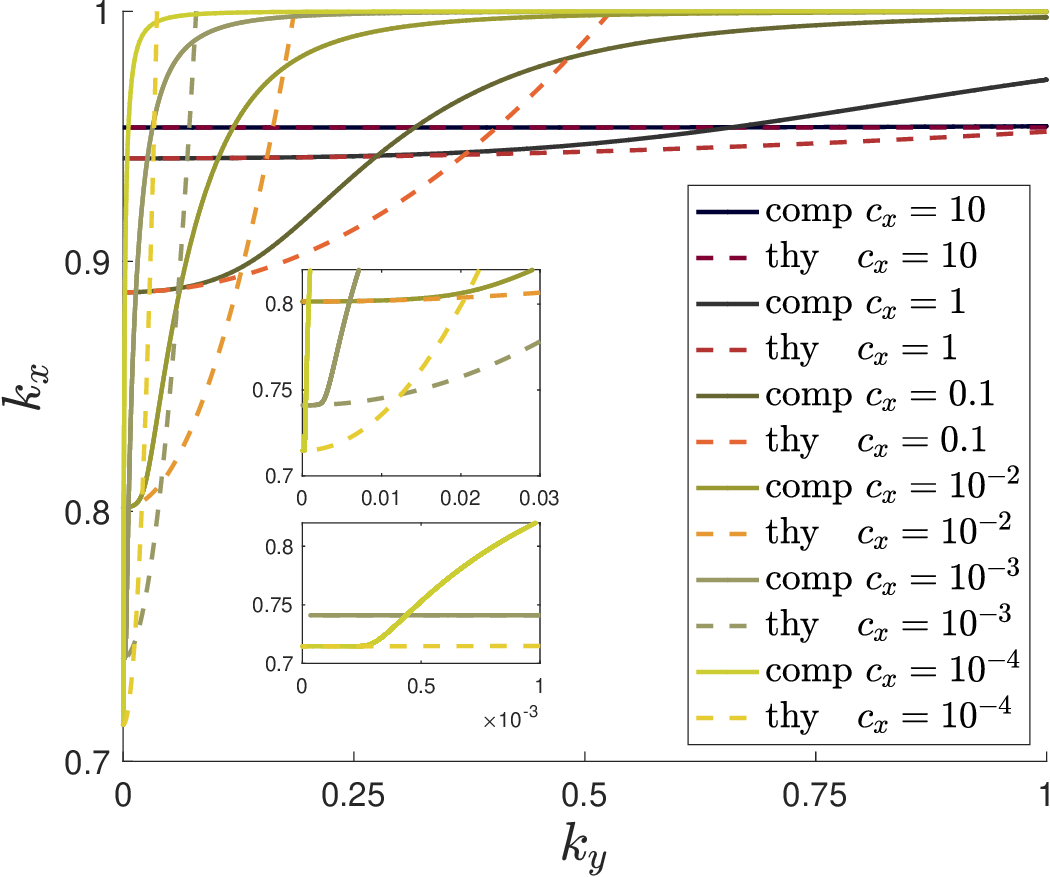}\hspace*{.01in}
 \includegraphics[width=0.195\textwidth]{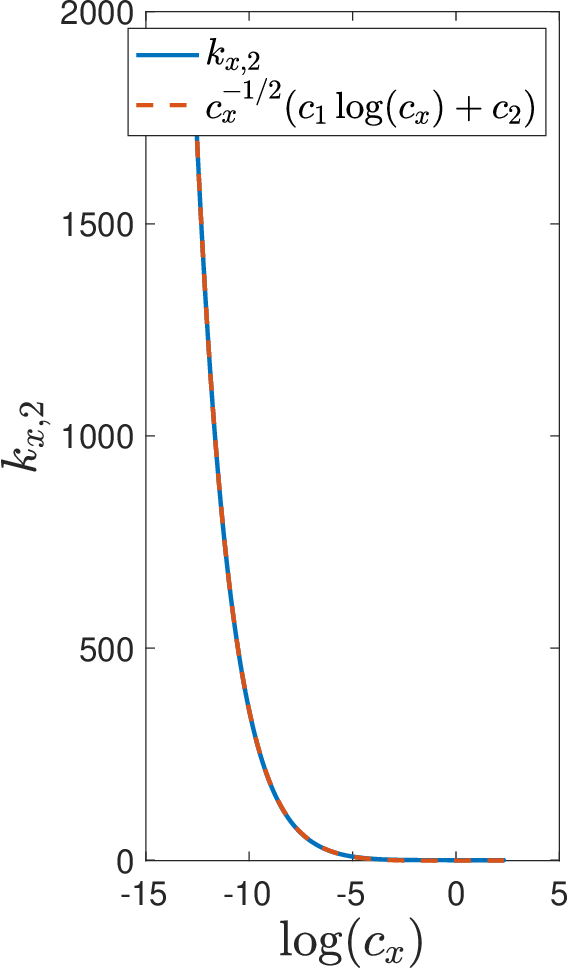}\hspace*{.01in}
  \includegraphics[width=0.385\textwidth]{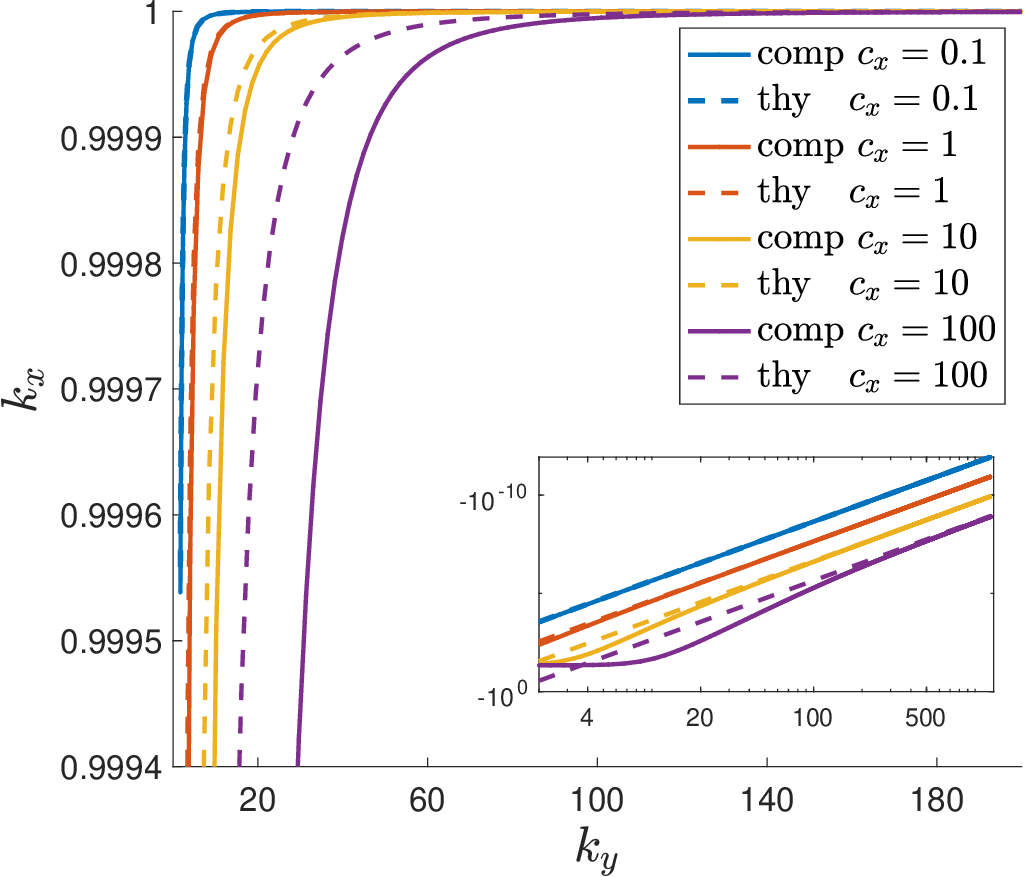}
 \caption{Left: Selected $k_x$ vs $k_y$ for $k_y$ small, compared with numerically computed quadratic approximation \eqref{e:kyto0}; note the good fit, albeit on increasingly small $k_y$-ranges as $c_x$ decreases. Center: Quadratic coefficient $k_{x,2}$ in  \eqref{e:kyto0} vs $c_x$ and comparison with best fit $c_1=-0.3422$, $c_2=-1.0439$ in \eqref{e:cx2asy}. Right: Selected $k_x$ vs $k_y$ for $k_y$ large and sample values of $c_x$, compared with \eqref{e:kx3}; inset log-log plot of $1-k_x$ vs $k_y$ confirming the good cubic approximation for moderate values of $c_x$.}\label{f:3}
\end{figure}

Expanding in inverse powers $\eps=1/k_y$, we find formally at orders $-1,0,1$, 
\begin{align}
 \rmO(-1):&\  |\partial_\zeta|\psi_0=0,\nonumber\\%label{e:oo-1}\\
 \rmO(0):& \ |\partial_\zeta|\psi_1+k_{x,0}-\frac{c_x}{2} \psi_0 -g(\psi_0+\zeta)=0,\nonumber\\%\label{e:oo0}\\
 \rmO(1):& \ |\partial_\zeta|\psi_2-\frac{1}{2}c_x\psi_1 +\frac{1}{8}|\partial_\zeta|^{-1}\left(c_x^2 + 4 c_xk_{x,0} \partial_\zeta\right)\psi_0+\left(k_{x,1}-g'(\psi_0+\zeta)\psi_1\right)=0.\nonumber%\label{e:oo1}
\end{align}
At $\rmO(-1)$,  we set $\psi_0=0$, which gives at $\rmO(0)$, 
\[k_{x,0}=\dashint g, \quad \psi_1=|\partial_\zeta|^{-1}(g-\dashint g).\]
Substituting the result into the equation at $\rmO(1)$ yields
\[
 |\partial_\zeta|\psi_2-\frac{1}{2}c_x\psi_1 +\left(k_{x,1}-g'(\zeta)\psi_1\right)=0,
\]
which upon averaging gives 
\[
 k_{x,1}=\dashint g'(\zeta)|\partial_\zeta|^{-1}(g(\zeta)-\dashint g)=0,
\]
which can be readily seen upon expanding $g$ in Fourier series, and 
\[
 \psi_2=|\partial_\zeta|^{-1}\left((g'+\frac{1}{1}c_x)\psi_1\right).
\]
Assuming that $g'$ is even, for instance $g=1+\kappa \sin(v)$, we see that $\psi_1$ and $\psi_2$ are both odd. 
At the next order, we find 
\[
 k_{x,2}=\dashint \left(
 (-\frac{1}{2}c_x-g'(\zeta))\psi_2 -4g''(\zeta)(\psi_1)^2
 \right),
\]
which vanishes when $g'$ is even. Continuing further the expansion, we find that the even part of $\psi_3$ is nonzero, 
\[
 \psi_{3,\mathrm{e}}=|\partial_\zeta|^{-3}\left(
 -\frac{1}{2}c_xk_{x,0} g'(\zeta)
 \right),
\]
and therefore 
\[
k_{x,3}=\dashint \psi_3g'\neq 0.
\]
In the specific case $g(v)=1+\kappa \sin(v)$, we find 
\begin{equation}\label{e:kx3}
 k_x=1 +k_{x,3} k_y^3 + \rmO(k_y^4),\qquad k_{x,3}=-\frac{1}{4}c_x\kappa^2;
 \end{equation}
see Figure \ref{f:3} for comparison with directly computed solutions. Note in particular that the asymptotics become steeper as $c_x$ increases, accommodating thus for the mismatch of limiting values, 
\[
 \dashint g = \lim_{k_y\to \infty}\lim_{c_x\to\infty} k_x \neq 
 \lim_{c_x\to\infty}\lim_{k_y\to \infty} k_x =\left(\dashint g^{-1}\right)^{-1};
\]
compare also the graphs in Figure \ref{f:all}.

\subsection{Qualitative summary and numerical explorations}\label{s:sum}

In the specific case of $g(v)=1+\kappa\sin(v)$, the asymptotics described above coincide well with numerical computations and predictions from the asymptotics give a good qualitative overall picture. 

\paragraph{Behavior for fixed $k_y$} Fixing $k_y$ small, we discuss the curve $k_x(c_x)$. By the integral identities above, $k_x(0)=1$ and $k_x'(0)<0$, while $k_x(\infty)<1$, monotonically increasing for $k_y<k_y^*$ and monotonically decreasing for $k_y>k_y^*$, $c_x\gg 1$. The asymptotics are therefore compatible with  globally monotonically decreasing $k_x(c_x)$ for $k_y>k_y^*$ and with $k_x(c_x)$ having a unique minimum for some finite $c_x(k_y)$ for $k_y<k_y^*$. This simple behavior with unique minimum or simple monotonicity is indeed what we observe numerically.

\paragraph{Behavior for fixed $c_x$} From the analysis above, we found $k_x(\infty)=1$ and $k_x$ monotonically increasing for large $k_y$ \eqref{e:kx3}. For $k_y=0$, the asymptotics and numerical analysis in \cite{beekie} predict $1-\kappa<k_x(0)<1$. The asymptotics with numerical evaluations of the relevant integrals predict that $k_x$ is monotonically increasing for $k_y\sim 0$, as well. Curves $k_x(k_y)$ computed numerically are in fact monotonically increasing on $k_y\geq 0$, albeit with a characteristic transition that we will discuss in the next section. 

\paragraph{Behavior as $k_y\to 0$} One notices that the limit of curves $k_x(c_x)$ as $k_y\to 0$ is not regular. In fact, at $k_y=0$, the results in \cite{beekie} show a monotone curve $k_x=1-\kappa+\rmO(\sqrt{c_x})$, and $k_x\in [1-\kappa,1+\kappa]$ for $c_x=0$. For $k_y>0$ curves $k_x(c_x)$ are non-monotone and appear to converge to this limiting set $(c_x,k_y)\in 0\times [1-\kappa,1+\kappa] \cup \{(c_x,k_x(c_x)),\ c_x>0\}$.  

\paragraph{Summary} Rephrasing our findings in terms of strain, measured through the deviation of $k_x$ from the equilibrium strain $k_x=1$, induced on stripes through forced growth at rate $c_x$ and imposed angle determined by $k_y$, we can summarize our findings as follows. 
\begin{enumerate}
 \item for small angles, $k_y\sim 0$, slow growth creates the largest residual strain in the stripes. For zero angles, $k_y=0$, the strain decreases with  increased growth rate, but for small angles the residual strain first increases with $c_x$ before faster growth reduces strain;
 \item for fixed growth rate, residual strain decreases with increasing angles;
 \item for larger angles, strain increases with growth rate. 
\end{enumerate}

The induced strain at $k_y=0$ can be understood as a non-adiabatic effect, proportional to $\kappa$ which measures the non-adiabaticity, that is, the size of terms that do not commute with the phase averaging symmetry $\varphi\mapsto \varphi+const$. Stripes are stretched maximally for small speeds, repeated stripe nucleation helps release stress with increased growth rate as described in \cite{beekie}. For small angles, an effect similar to zero angle can be observed, with the caveat that for very small speeds, the gliding of a localized boundary defect along the growth interface can mediate the growth process with little residual stress. Increasing the rate of growth increases the glide speed of the defect and thereby residual strain. Yet stronger growth leads to a phase transition in the nature of the boundary defect that leads to delocalization and decreased strain.  

\begin{figure}[h]
 \centering

 \includegraphics[width=0.33\textwidth]{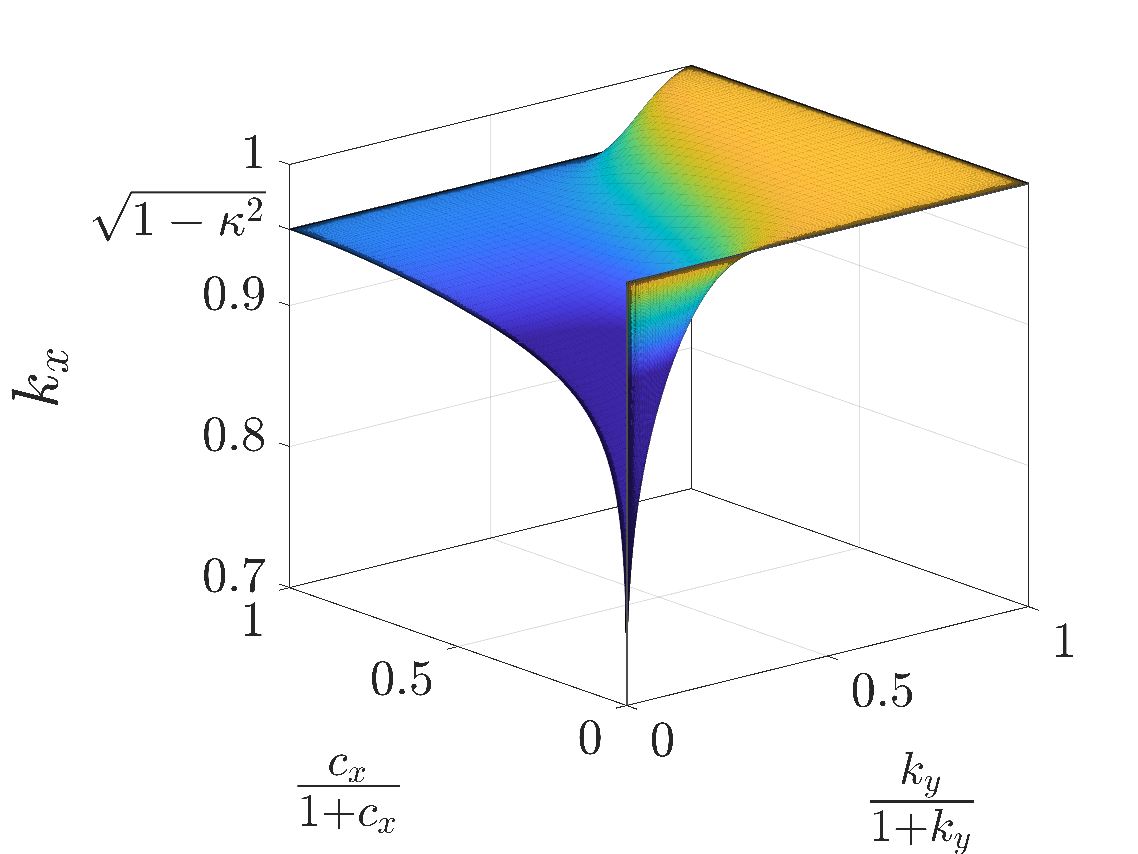}%
 \includegraphics[width=0.33\textwidth]{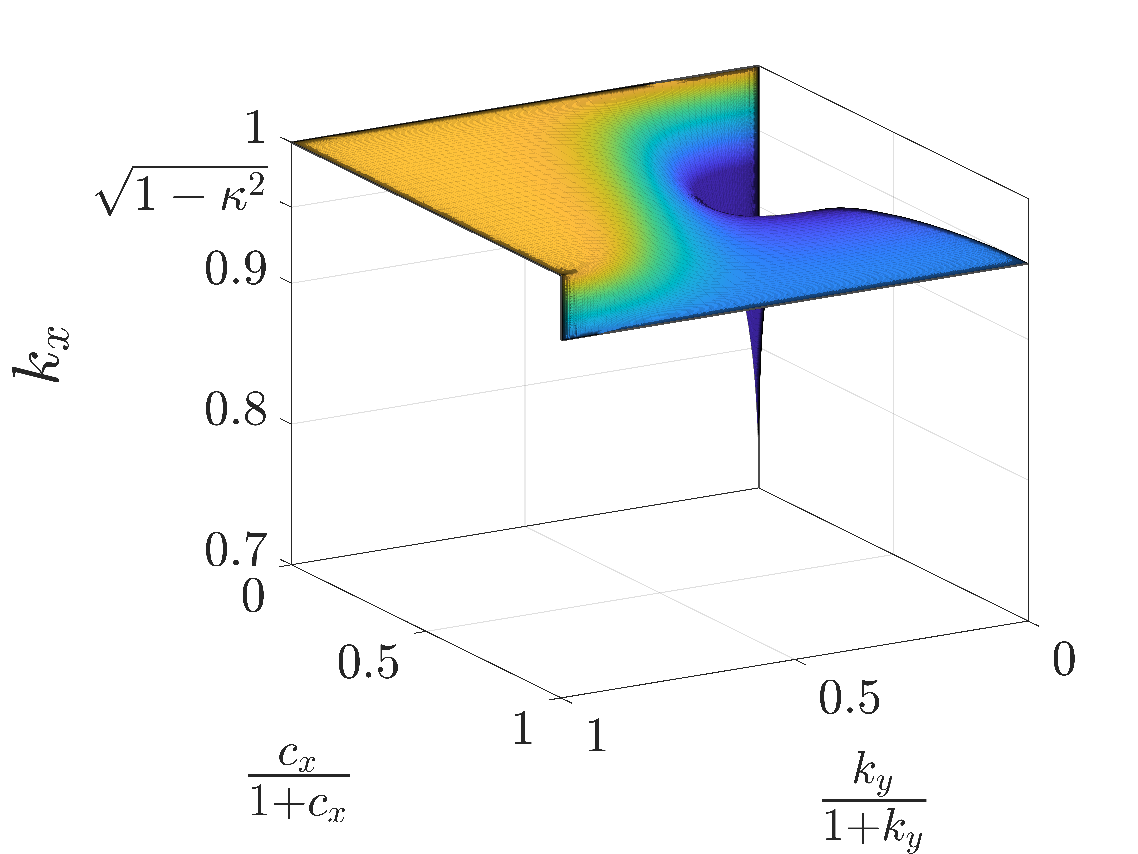}% 
 \includegraphics[width=0.33\textwidth]{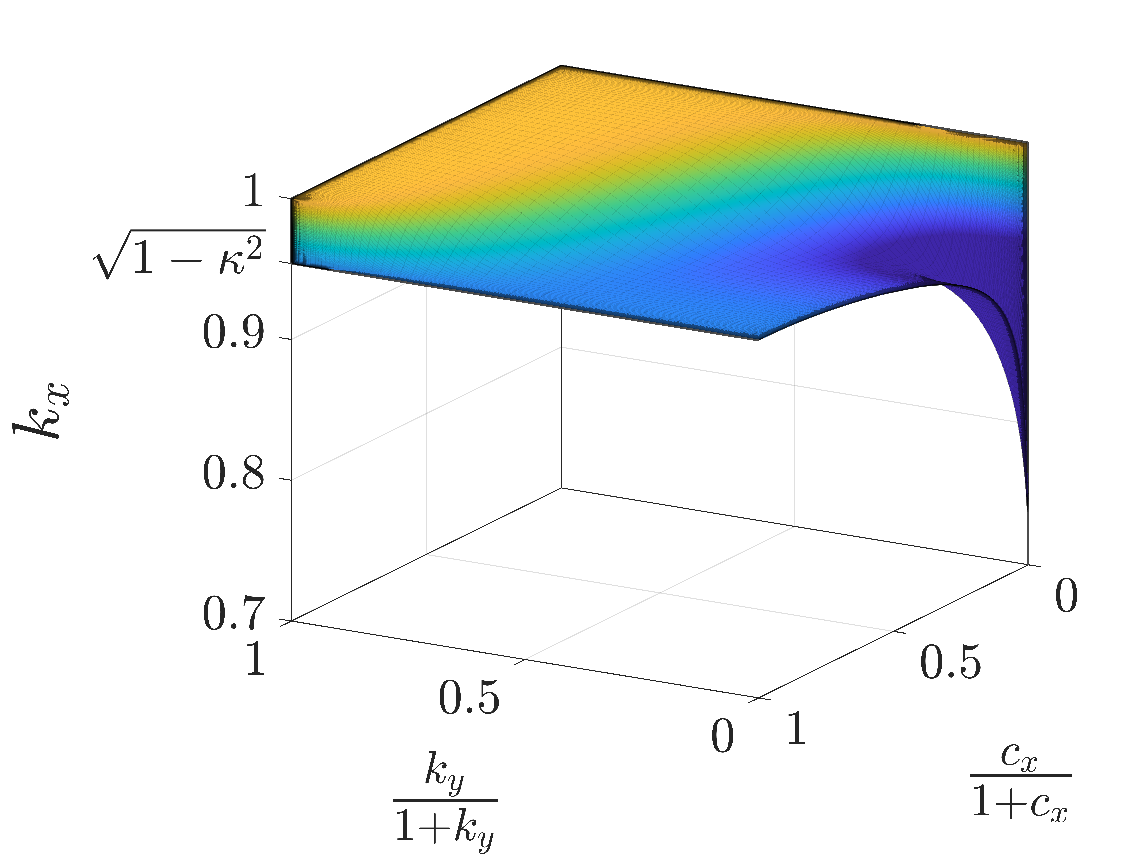}% 
 \caption{Surface $k_x$ as a function of $k_y$ and $c_x$. Plots use $k_y/(1+k_y)$ and $c_x/(1+c_x)$ as coordinates to include the limits $c_x=\infty$ and $k_y=\infty$ at 1; see also \texttt{mod\_space\_all.mp4} in the supplementary materials. }\label{f:all}
\end{figure}

Increasing the angle through $k_y$ reduces the non-adiabaticity, up to the point where stripes perpendicular to the boundary can grow without deformation at the interface, $k_y=1$, not creating any strain. Figure \ref{f:all} shows the surface $k_x(k_y,c_x)$ from different angles, exhibiting the singularities that occur in the compactification at the boundaries $c_x,k_y\in \{0,\infty\}$.

\section{Asymptotics near the origin}\label{s:origin}
The strain in a large region of parameter space is simply monotone and fairly simple asymptotics explain the behavior. The most intriguing, non-monotone dynamics occur in a vicinity of $c_x=k_y=0$. In this regime, profiles $\varphi$ converge to step-like functions in $\zeta$; see Figure \ref{f:profiles}. An inner expansion of the layer-type solution reveals an interesting transition that sheds light on the asymptotics in this region. 

We scale in \eqref{e:pdp1}--\eqref{e:pdp4} for an inner expansion at the heteroclinic  $k_y=\tilde{k}_y\eps$, $c_x=\eps$ and  $\partial_\zeta=\eps \partial_z$, and obtain, expanding the Fourier symbol $\mathcal{D}_+$, at leading order
\begin{equation}\label{e:inner}
  \mathcal{D}\psi =g(\psi) - k_x, \quad y\in\R,\qquad  \psi(-\infty)+2\pi=\psi(+\infty)=\psi_*, \quad \mathcal{D}=\sqrt{-\tilde{k}_y^2\partial_{zz}+k_x\partial_z}, 
\end{equation}
where $\mathcal{D}$ now is defined as a Fourier multiplier for functions on the real line rather than periodic functions. This equation does  have a local interpretation as a traveling-wave solution $\psi=\psi(\tilde{k}_y y + k_x t,x)$ to the heat equation with nonlinear boundary flux,
\[
 \psi_t=\Delta \psi,\ x<0,y\in\R,\qquad \psi_x=g(\psi)-k_x,\ x=0,y\in\R.
\]
Such traveling waves have been studied in \cite{cabre}, establishing in particular existence and monotonicity properties for solutions $\psi(y-ct)$, with $c=c(k_x)$ for $|k_x-1|<\kappa$ when $g(\psi)=1+\kappa \sin(\psi)$. Rescaling $y=z/k_y$ shows that these traveling solutions give solutions to \eqref{e:inner} whenever 
\begin{equation}\label{e:kykx}
 k_y=\frac{k_x}{c(k_x)}.
\end{equation}
Moreover, monotonicity of $c$ in $k_x$ from \cite{cabre} implies that $k_y$ is monotonically increasing as a function of $k_x$ with minimum $k_y^*$, such that we can rewrite \eqref{e:kykx} as 
\begin{equation}\label{e:kxky}
 k_x=k_x^\mathrm{f}(k_y), \qquad \text{for} \ k_y>k_y^*. 
\end{equation}
For $0<k_y<k_y^*$, we conjecture the existence of heteroclinic solutions with $k_x=\min g(\varphi)$, asymptotic to $\mathrm{argmin}\, g(\varphi)$ and $\mathrm{argmin}\, g(\varphi)+2\pi$. In particular, the selected $k_x$ is constant at leading order. 

Below, we provide numerical evidence for our predictions. 

\paragraph{Computing heteroclinic orbits in \eqref{e:inner}}

We focus on the specific case $g(\psi)=1+\kappa\sin\psi$. 
In order to solve \eqref{e:inner}, we rely on Fourier transform. We therefore write $\psi=\psi_\mathrm{s}+\tilde{\psi}$ with  $\psi_\mathrm{s}(z)=\psi_*+2\arctan(z)$.  The profile $\psi_\mathrm{s}(z)$ accounts for the heteroclinic structure such that $\tilde{\psi}$ can be chosen to be periodic. The asymptotic state is (necessarily) chosen such that $g(\psi_*)=k_x$, $g'(\psi_*)\geq 0$. The choice of  $\arctan(z)$ is motivated by the fact that the action of the integral operator is explicit, 
\[
 \mathcal{R}(z;k_x,\tilde{k}_y):=\mathcal{D}\psi_\mathrm{s}(z)=\frac{2\sqrt{\pi}\tilde{k}_y}{1 + z^2}
 \Re\left((1+\rmi)\mathrm{U}\left(-\frac{1}{2},0,\frac{k_x(-\rmi +z)}{\tilde{k}_y^2}\right)\right),
\]
where $\mathrm{U}$ is the confluent hypergeometric Kummer-U function. We then solve 
\[
 \mathcal{D}(k_x,\tilde{k}_y)\tilde{\psi}+\mathcal{R}(k_x,\tilde{k}_y)-g(\psi_\mathrm{s}+\tilde{\psi})+k_x=0,
\]
with periodic boundary conditions on a large domain $|z|\leq L$ together with a phase condition 
$\int\tilde{\psi}(z)\rme^{-z^2}\rmd z=0$ and with $k_x$ as a Lagrange multiplier using a Newton method and secant continuation in $k_y$. The spectral discretization gives accuracy of $10^{-6}$ for moderate effective discretization sizes of $0.1$. Solutions decay however only weakly with $z^{-1/2}$, $z\to -\infty$, and $z^{-3/2}$ for $z\to+\infty$. We found accuracy of $10^{-6}$ for domain sizes $L\sim 10^6$ using $N=2^{24}\sim 10^7$ Fourier modes. The code was implemented in \texttt{matlab} and ran on an Nvidia GV100 graphics card allowing for fast evaluation of the large discrete Fourier transforms. The Kummer-U function was evaluated and tabulated in \texttt{mathematica} and interpolated in \texttt{matlab}, since direct evaluation in \texttt{matlab} is slow. 

 \begin{figure}[h]\centering\includegraphics[width=0.54\textwidth]{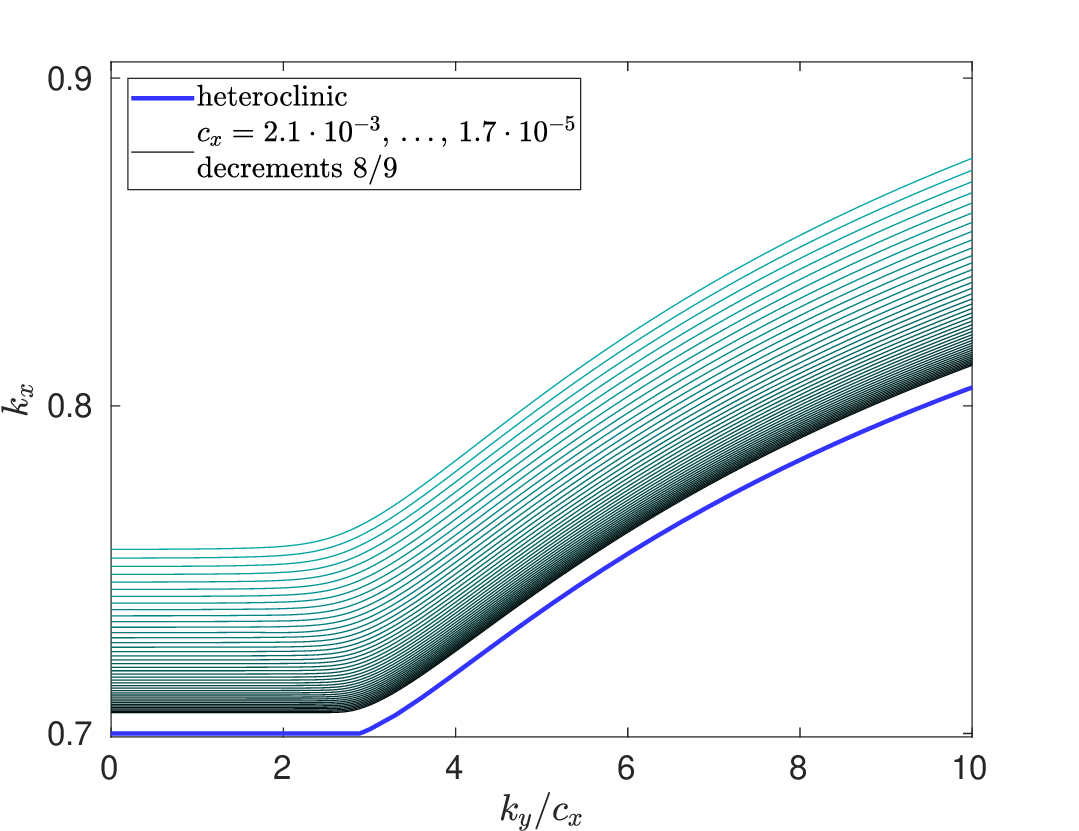}\hspace*{-0.2in}
 \includegraphics[trim=0.25in 0 0 0, clip,width=0.49\textwidth]{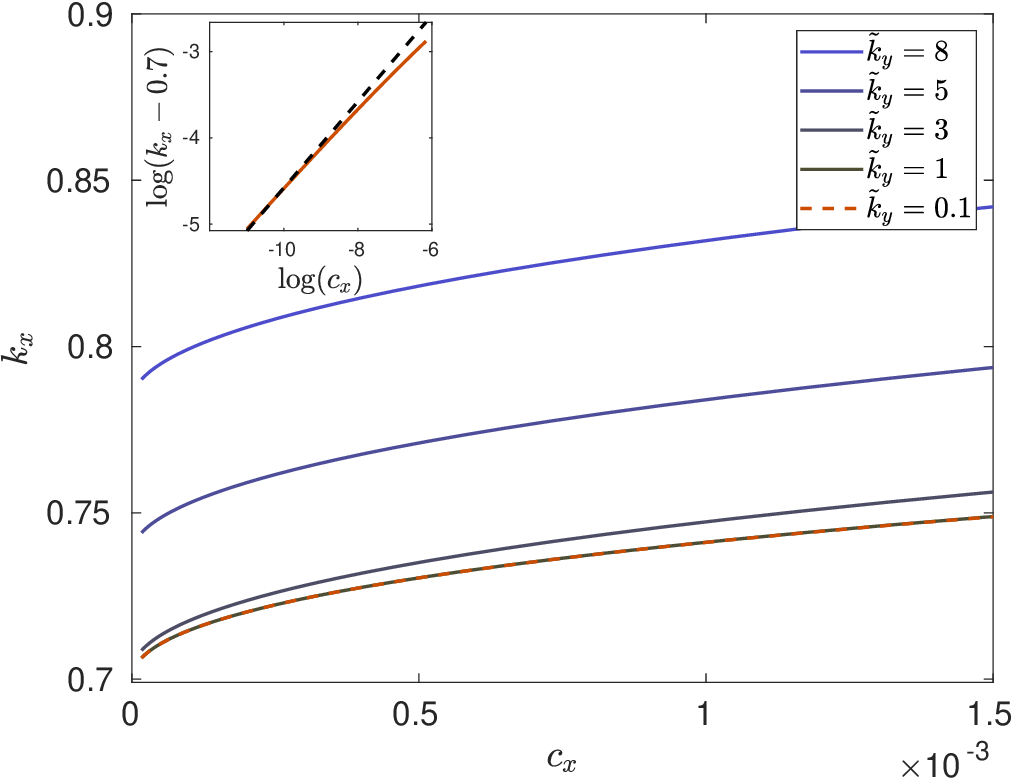}
 \caption{Left: selected $k_x$ in \eqref{e:bi} plotted against $k_y/c_x=\tilde{k}_y$, for $c_x$ decreasing geometrically by factors $8/9$ and with the selected $k_y$ from the heteroclinic continuation for comparison.  Right: section through the diagram on the left, plotting $k_x$ as a function of $c_x$ for fixed $\tilde{k}_y$, showing in particular that the values are almost independent of $\tilde{k}_y<2$, below the heteroclinic bifurcation. For such small values, $k_x \sim 0.7+1.52\sqrt{c_x}$ in good agreement with \cite{beekie}.}\label{f:4}
\end{figure}
\begin{figure}[h]
\begin{minipage}{0.43\textwidth} \centering\includegraphics[width=\textwidth]{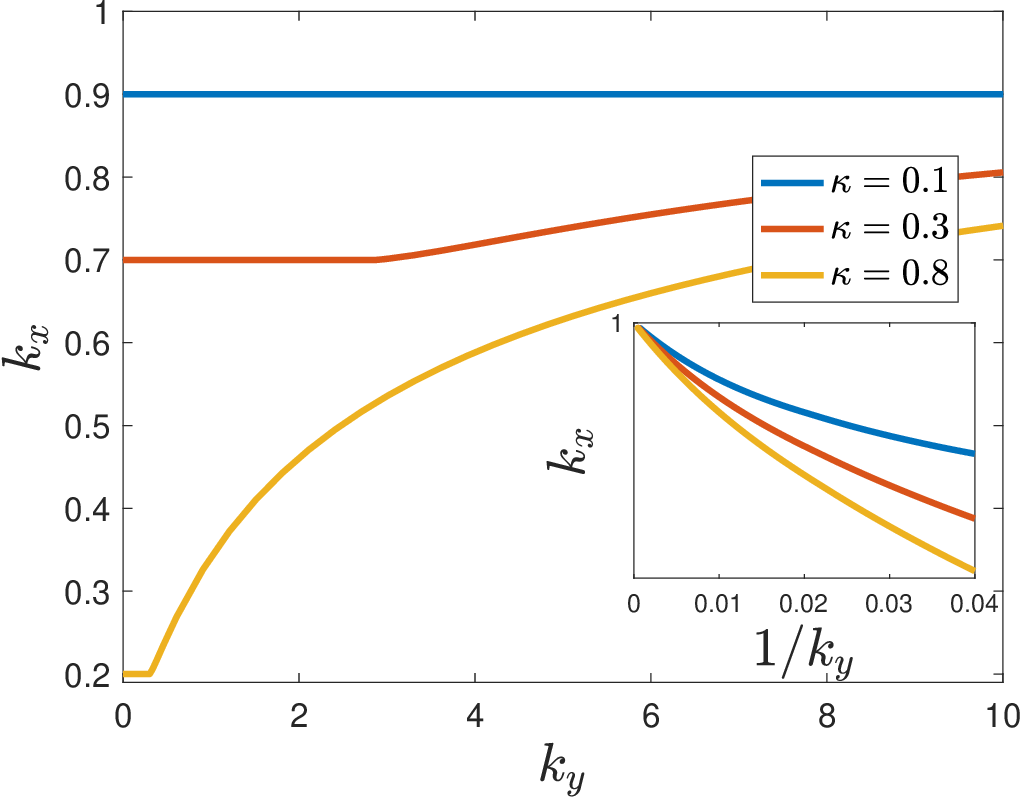}\\
 \includegraphics[width=\textwidth]{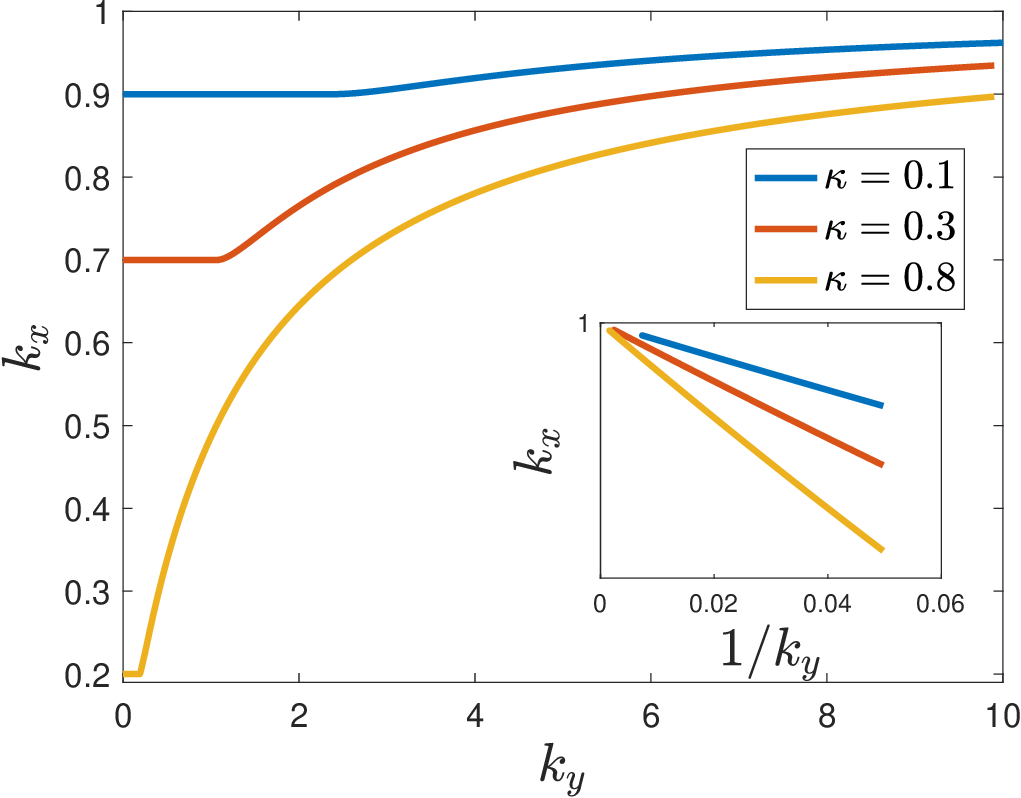}
 \end{minipage}\hfill
\begin{minipage}{0.45\textwidth} %\centering 
\includegraphics[width=\textwidth]{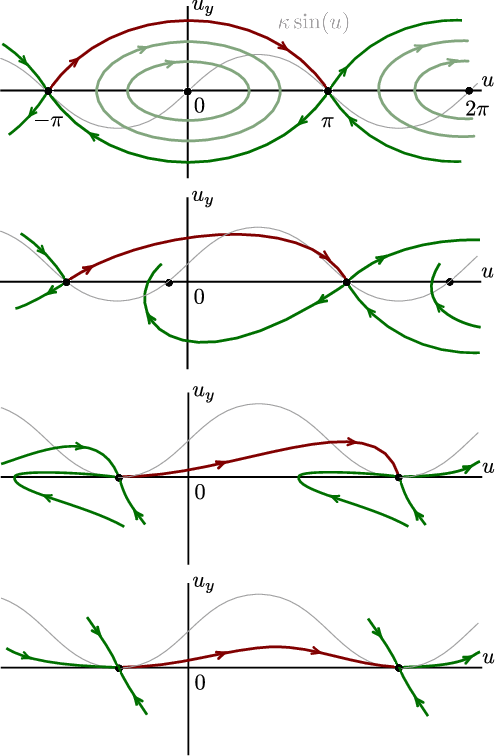}
\end{minipage} %  \includegraphics[width=0.18\textwidth]{figures/kx2}\hspace*{.1in}
 \caption{ Selected wavenumbers in nonlocal \eqref{e:inner} and local \eqref{e:loc} problems (left top and bottom, resp.) for several values of $\kappa$. Both show the distinct transition to a flat regime for small $k_y$, where the nature of the heteroclinic changes and prevents further increase of the deviation from equilibrium strain $k_x=1$.  The transition in the local case can be understood as a heteroclinic flip bifurcation with phase portraits depicted on the right, with the flat piece of the $k_x--k_y$ graph corresponding to the saddle-node heteroclinic at the bottom and the transition occurring at the heteroclinic flip bifurcation. }\label{f:het}
\end{figure}

Results from the computation of heteroclinic orbits are shown in Figure \ref{f:het}, left upper panel, showing a characteristic transition from increasing values $k_x(\tilde{k}_y)$ for moderate $\tilde{k}_y$ to constant $k_x$ for small $\tilde{k}_y$. At the transition value, the heteroclinic orbit delocalizes, the amplitude of $\psi_y$ decreases. In the limit $k_y\to\infty$, we find the ``Hamiltonian'' picture, with $k_x=1$. 

The computed values of $k_x$ compare well with the selected values in the selection problem periodic in $y$, as shown in Figure \ref{f:4}. Selected wavenumbers $k_x$ as function of the scaled wavenumber $k_y/c_x$, computed for fixed values of  $c_x\ll1 $ through continuation in $k_y\to 0$,  converge to the limiting curve given by the heteroclinic orbit. 

The nonlocal  problem is related to the Weertman equation that is used to describe the glide motion of dislocations; see \cite{josien} and references therein. In fact, the nonlocal Weertman equation can be obtained by replacing our nonlinear fluxes by a dynamic (Wentzel) boundary conditions,
\[
\varphi_t=\Delta \varphi,\ x<0;\qquad  \varphi_t=-\varphi_x+g(\varphi),\ x=0.
\]
Our numerical methods in fact resemble the approach taken in \cite{lebris}, although pseudo-differential operators are more difficult in our case and the emphasis in \cite{josien} is on the time-dependent initial-value problem. 
We conclude this analysis with a heuristic explanation of the transition from a sharply localized defect selecting strains $k_x$ to a delocalized heteroclinic selecting minimal values of $k_x$, through analogy to a local differential equation.

\begin{figure}[h]
 \centering
 \includegraphics[width=1\textwidth]{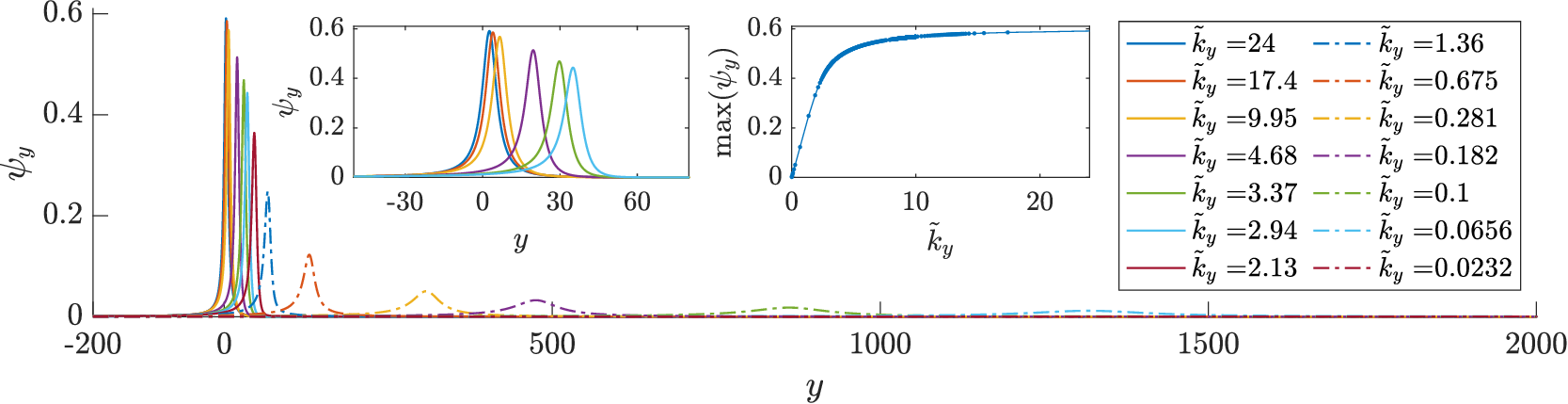}\\
  \includegraphics[width=1\textwidth]{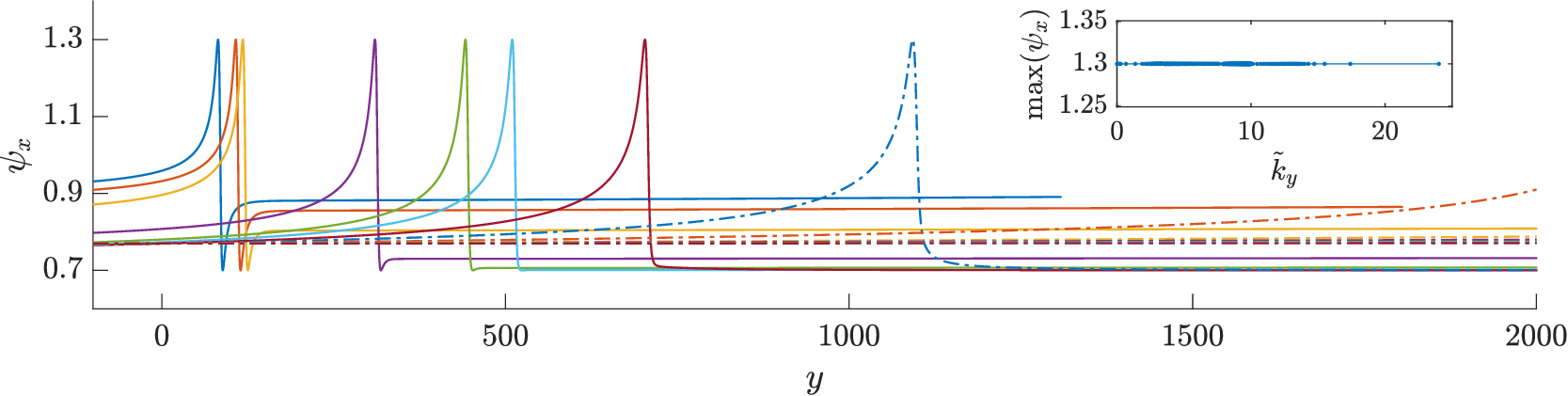}
 
 \caption{Profiles of derivatives $\partial_y\varphi$ (top) and $\partial_x\psi$ in unscaled variables $y$ for values of $\tilde{k}_y=k_y/c_x$, $c_x=10^{-4}$, passing the heteroclinic bifurcation. Profiles are roughly constant for large $\tilde{k}_y$ (left inset) but rapidly delocalize past the heteroclinic transition with long tails to the left of the peak; amplitude of profiles rapidly decreases past heteroclinic bifurcation (right inset). Normal derivatives also delocalize but always peak at minimal and maximal strain. }\label{f:profile_transition}
\end{figure}

\paragraph{Comparison with local heteroclinic bifurcations}
A qualitatively equivalent picture emerges when the nonlocal pseudo-differential operator $\mathcal{D}$ is replaced by a local operator $\mathcal{D}_\mathrm{loc}=-\tilde{k}_y^2\partial_{zz}+k_x\partial_z$. In this case, elementary phase plane analysis establishes the existence of heteroclinic orbits to $\mathcal{D}_\mathrm{loc} \psi= 1+\kappa \sin(\psi)-k_x$. Rescaling $k_y\partial_z=\partial_y$, we find the traveling-wave equation to the (asymmetric) parabolic Sine-Gordon equation,
\begin{equation}\label{e:loc}
 u_{yy}+c u_y=1+\kappa \sin(u)-k_x,\qquad  c=k_x/k_y.
\end{equation}
For $k_x=1$ we have $c=0$ and a heteroclinic between $u=0$ and $u=2\pi$. The heteroclinic is transversely unfolded in the parameter $c$ and we can in fact continue the heteroclinic with $c=c(k_x)$ monotonically increasing as $k_x$ is decreasing, until $k_x=1-\kappa$.  For $c\gg 1$, we find at leading order, after a reduction to a slow manifold, 
\[
c u_y=1+\kappa \sin(u)-k_x,
\]
which possesses heteroclinic orbits for $k_x=1-\kappa$, connecting the saddle-node equilibria $u=-\pi/2\mod 2\pi$. These heteroclinics between saddle-node equilibria are robust up to a heteroclinic codimension-two bifurcation \cite{chowlin,bellay}. The associated phase-portraits in the $u-u_y$--plane are shown in Figure \ref{f:profile_transition} and can be easily confirmed using elementary phase-plane analysis and monotonicity in $c$. 

%%%%
%%%%
%%%%
\section{Comparison with an anisotropic Swift-Hohenberg equation}\label{s:anSH}
Returning to the motivation by striped patterns, we now study the formation of striped patterns in a directionally quenched Swift-Hohenberg equation.
The phase-diffusion approximation with nonlinear boundary fluxes given by the strain-displacement relation was shown to be a correct approximation in the case $c_x=0$ in \cite{scheelweinburd2017}, for $y$-independent patterns.  Considering patterns in two spatial dimensions, one notices that patterns selected for $c_x\ll 1$ and $k_y\ll 1$ have wavenumber $k<1$ and are zigzag unstable; see again, for instance, \cite{scheelweinburd2017}. As a consequence, a phase-diffusion approximation for dynamics of these patterns would yield a negative effective diffusion coefficient in the direction along stripes and higher-order corrections as in the Cross-Newell equation are necessary to fully capture dynamics; see for instance \cite{pismen}.

We therefore focus on the  quenched anisotropic Swift-Hohenberg equation,
\begin{equation}\label{e:aSH}
u_t = -(1+\Delta_{x,y})^2 u + \beta \partial_{yy} u + \mu u - u^3,
\end{equation}used in \cite{boyer,https://doi.org/10.1029/2010GL046091,handwerk,PhysRevE.69.066213,fl,MIKHAILOV200679,MUNOZGARCIA20141,PhysRevE.73.036117,PhysRevE.65.046219} to describe nematic liquid crystals, electroconvection, ion bombardment, surface catalysis, or vegetation patterns; see also \cite{haragus2012dislocations} for an analysis of dislocations in this model. 
For $\beta>0$, the anisotropic term suppresses the zig-zag instability in stripes with wavenumbers $k\lesssim1$. For sufficiently large $\beta$ all wavenumbers within the strain-displacement relation, $k\in(k_\mathrm{min},k_\mathrm{max})$, with $k_\mathrm{max} = \max g(\phi)$, are stabilized. In the following, we first derive a phase-diffusion approximation and nonlinear fluxes in the form studied in this paper from the anisotropic Swift-Hohenberg equation, and then describe a numerical approach to computing striped patterns created in directional quenching, with the goal of comparing the numerical results to the quantitative predictions from the phase-diffusion approximation. Throughout, we focus on the regime $0<c_x,k_y\ll 1$ and use a quenched parameter of the form $\mu=-\mu_0\tanh((x-c_x t)/\delta)$ with $\delta=0.5$.

\paragraph{Derivation of phase diffusion in anisotropic Swift-Hohenberg}
Focusing on nearly parallel stripes with constant parameter $\mu$, we use the parabolic scaling $\mu = \epsilon^2,  x = \epsilon\tilde x ,\,  y = \epsilon \tilde y,\, t = \epsilon^2 \tilde t$, and substitute the ansatz $u(x,y,t) = \eps A\left(\tilde x,\tilde y, \tilde t\right) e^{\rmi  x} + \mathrm{c.c.}$ into \eqref{e:aSH} to obtain, at leading order, an anisotropic Ginzburg-Landau equation
 \begin{equation}
A_{\tilde t} = 4 A_{\tilde x\tilde x} + \beta A_{\tilde y\tilde y} + A - 3 A|A|^2. 
\end{equation}
Introducing polar coordinates $A=R\rme^{\rmi\tilde{\phi}}$ and expanding near $R=1/\sqrt{3}$, $\tilde{\phi}=0$, one finds an exponentially damped equation for $R$ and an anisotropic diffusion equation for $\tilde{\phi}$, 
\begin{equation}\label{e:pds}
 \tilde{\phi}_{\tilde t}=4 \tilde{\phi}_{\tilde x\tilde x} + \beta \tilde{\phi}_{\tilde y\tilde y} .
\end{equation}
% Next,  inserting $A\left(\tilde x, \tilde y,\tilde t\right) = \left(1+\delta r\right)\re^{\ri \tilde \phi}$, with $r = r\left(\tilde x,\tilde y,\tilde t\right)$, $\tilde \phi=\tilde\phi\left(\tilde x,\tilde y,\tilde t\right)$, and $0<\delta\ll1$,  one obtains at leading-order a linear anisotropic phase diffusion equation. 
Note that this equation is again invariant under the parabolic scaling such that we may consider \eqref{e:pds} in the original coordinates $t,x,y$ to describe patterns in \eqref{e:aSH}. 

We next turn to the effect of the spatial quenching. At the order of the Ginzburg-Landau equation, one does not capture the non-adiabatic effects of the parameter jump. We use the expression for the strain-displacement relation from \cite{scheelweinburd2017} for the strain-displacement relation in the one-dimensional case, unaffected by the anisotropic term, $\tilde \phi_{x} = g_{\mathrm{SH}}(\tilde \phi) := 1 + \frac{\mu_0}{16}\sin 2\tilde\phi + \rmO(\mu_0^{3/2})$. The symmetry $\tilde{\phi}\mapsto \tilde{\phi}+\pi$ is present at higher orders, as well, and caused by the $u\mapsto -u$ symmetry in the nonlinearity and the ensuing symmetry $u_\mathrm{per}(\xi)\mapsto -u_\mathrm{per}(\xi+\pi)$ of periodic patterns. We use the same boundary condition for two-dimensional patterns, neglecting in particular dependence of $g_\mathrm{SH}$ on $\tilde{\phi}_y$, and also dependence on $c_x$, which gives the two-dimensional system
\begin{equation}\label{e:shpd}
\tilde \phi_{ t} = 4\tilde \phi_{ x  x} + \beta \tilde\phi_{ y  y} +  \tilde{c}_x\tilde \phi,\quad  x < 0,\,y\in\R ,\quad\qquad \tilde \phi_{ x} = g_{\mathrm{SH}}(\tilde \phi),\quad  x = 0,\,y\in\R. \end{equation}
With the additional scaling 
$
\phi = 2\tilde\phi,\quad x = \tilde x,\quad y = \tilde y,\quad  c_x = 8 \tilde c_x,\quad t = 16\tilde t,
$
we then  obtain the phase-diffusion equation \eqref{e:pd} with strain-displacement relation $\phi_x = g_{\mathrm{SH}}(\phi/2)$ at $x = 0$. We remark that by setting $\kappa = \mu_0/16$, $g_\mathrm{SH}$ agrees to leading order with the relation $\phi_x = g(\phi)$ employed in previous sections.  Through these scalings, we can compare the heteroclinic prediction of Section \ref{s:origin} with moduli curves of quenched patterned solutions $u(\tilde x,\tilde y,t)= u(k_x(\tilde x - \tilde c_x t), k_y(\tilde y - c_y \tilde t))$ of the full equation \eqref{e:aSH}. In our comparisons below, we use a value for $\kappa$ slightly different from $\mu_0/16$, computed directly from the one-dimensional Swift-Hohenberg equation as described in \cite{morrissey,scheelweinburd2017}, accounting for both error terms $\rmO(\mu_0^{3/2})$ and corrections due to the fact that we use a smoothed version of the step function for the spatially dependent parameter $\mu$. 

\paragraph{Oblique stripe formation in the full Swift-Hohenberg equation}
The formation of striped patterns is described by traveling-wave solutions \cite{gs4,zigzag} with speed vector $(c_x,c_y)$, again requiring $c_y = k_x \tilde c_x/ k_y$,
\begin{align}
0&= -(1+k_x^2\partial_{\xi}^2 + k_y^2 \partial_{\zeta}^2)^2 u + \beta k_y^2 \partial_{\zeta}^2 u+  \mu u - u^3 + \tilde c_x k_x(\partial_\xi + \partial_\zeta)u,\,\quad \xi<0, \zeta\in \R,\label{e:aSH2}\\
0&= u(\xi,\zeta+2\pi) -u(\xi,\zeta),\,\qquad \xi\leq 0, \zeta\in \R,\label{e:aSH3}\\
0&= \lim_{\xi\rightarrow\infty} u(\xi,\zeta) ,\qquad 0= \lim_{\xi\rightarrow-\infty} |u(\xi,\zeta) - u_\mathrm{per}\left(\xi+\zeta;k_x,k_y\right)|,\qquad \zeta\in \R.\label{e:aSH4}
\end{align} 
We numerically solve \eqref{e:aSH2} - \eqref{e:aSH4} using a farfield-core approach similar to \cite{lloydscheel,zigzag}, which decomposes $u = w + \chi u_\mathrm{per}\left(k_x,k_y\right)$, where $w$ is localized near the quenching interface, and $\chi$ is a cutoff function supported in the $\xi$-farfield. Here, we solve for $w$ and $k_x$ with parameter $k_y$, using a spectral discretization in both $\xi$ and $\zeta$ so that functions can be evaluated with the fast Fourier transform. Each Newton step of the pseudo-arclength continuation algorithm was once again performed using \texttt{gmres} to solve the associated linear problem. The nonlinear system was conjugated with exponentially localized weights and pre-conditioned with the principal symbol of the linear equation. Discretization and domain size were controlled adaptively ensuring both small tails at the end of the (periodic) domain and small amplitudes in highest Fourier modes. Typical domain sizes near the origin were $x\in (-800,800)$ with $8192\times4096$ Fourier modes in $(\xi,y)$.   Code was again implemented in \texttt{matlab} with computations carried out using an Nvidia GV100 GPU. Further details of this numerical approach are left for a companion work. For values of $\tilde{c}_x$ and $k_y$ smaller than the ones shown, \texttt{gmres} would usually not converge due to constraints on the number of inner iterations caused by limited memory. 

\begin{figure}[h]\centering
         \includegraphics[width=0.49\textwidth]{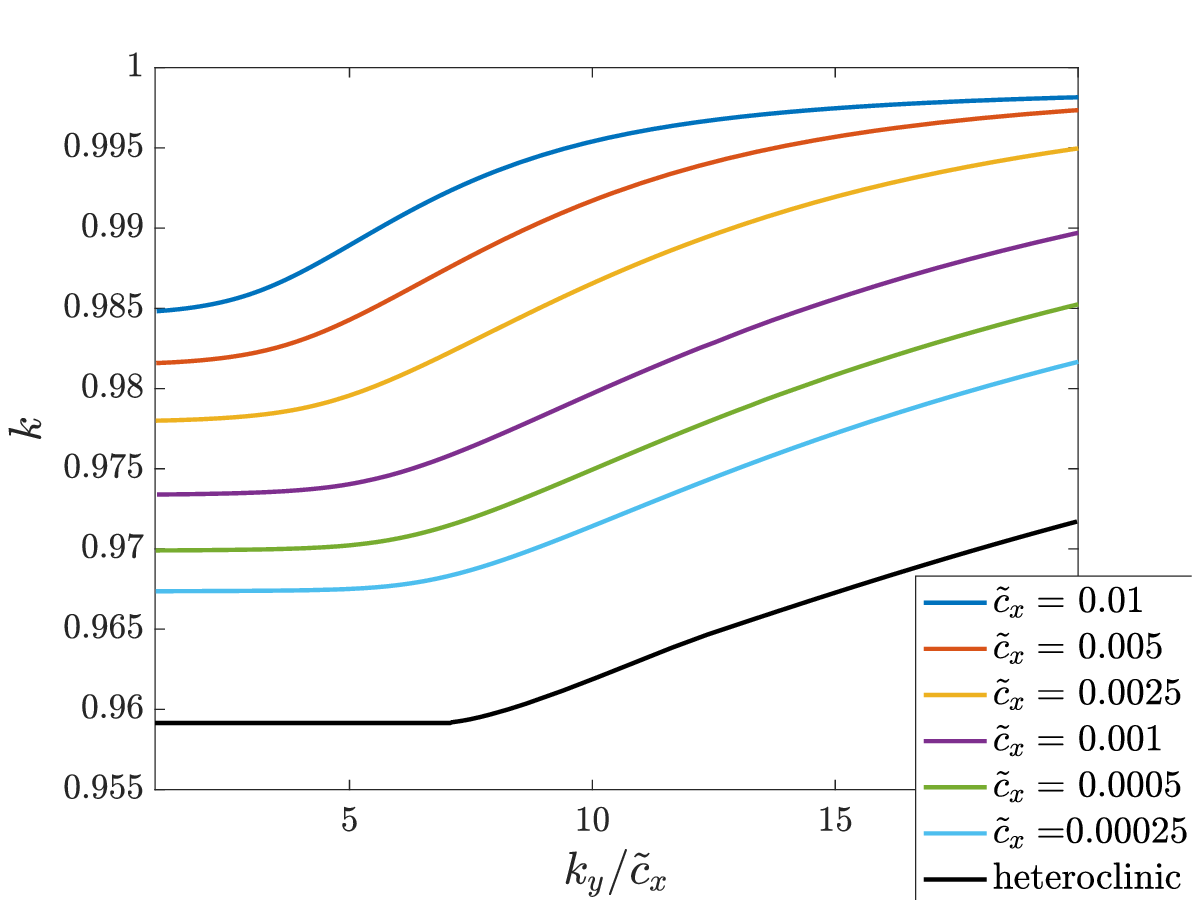}\hfill%width was .365
    \includegraphics[width=0.49\textwidth]{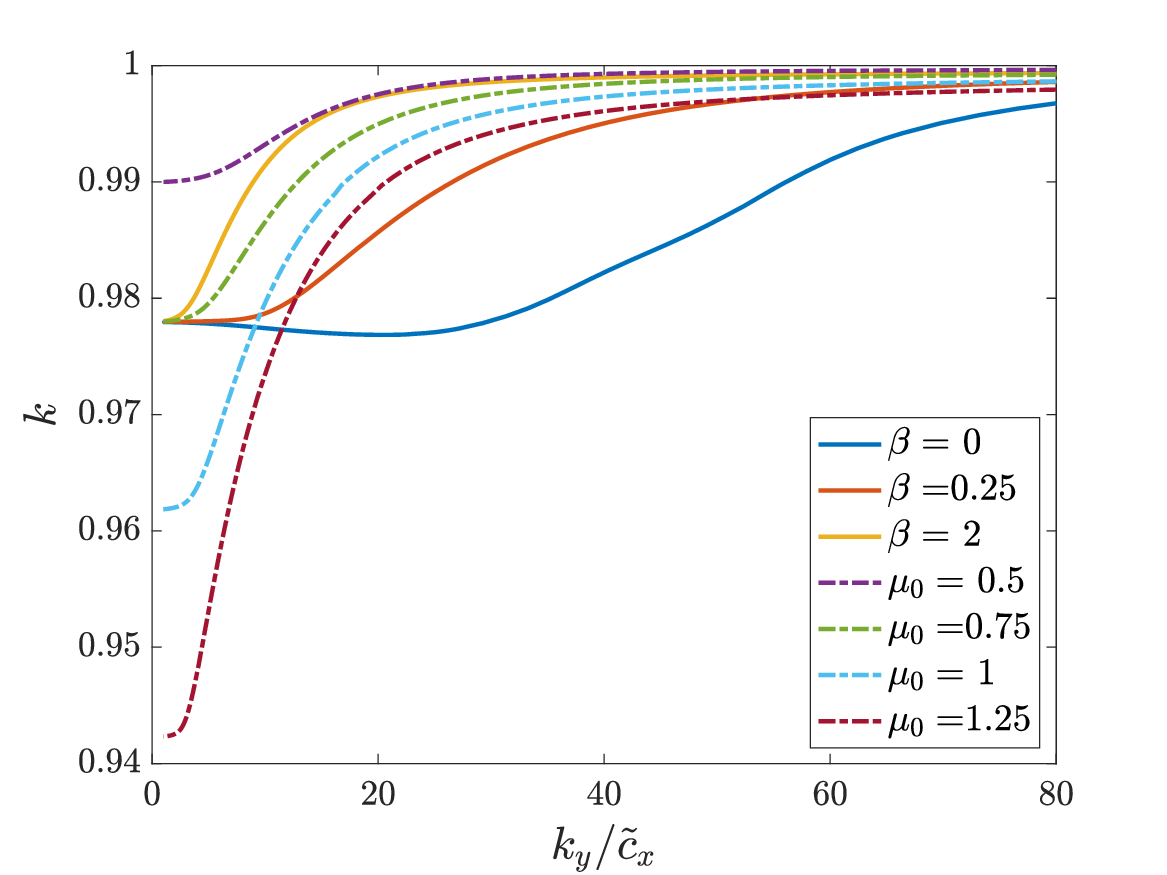}
     \caption{Wavenumber selection curves for anisotropic Swift-Hohenberg \eqref{e:aSH2}--\eqref{e:aSH4} for $k_y/\tilde c_x \sim 0$ with $\tilde c_x$ fixed. Left: comparison for  a range of $\tilde c_x$ values with the heteroclinic curve (black) of \S \ref{s:origin}; here, $\beta=1$ and $\mu_0 = 3/4$ so that $\kappa = \mu_0/16 = 3/64$. The heteroclinic curve (black) is obtained using numerically derived strain-displacement relation to account for higher-order corrections in $\mu_0$. Right: plot of selected wavenumber $k$ for $k_y/\tilde c_x \sim 0$ for a range of $\beta$ values with $\mu_0 = 3/4$ fixed (solid) and range of $\mu_0$ values with $\beta = 1$ fixed (dot-dashed), $\tilde c_x = 0.0025$.
}\label{f:aSH}
\end{figure}

\begin{figure}[h]
    \includegraphics[trim=0 0 0 0,clip,width=0.19\textwidth]{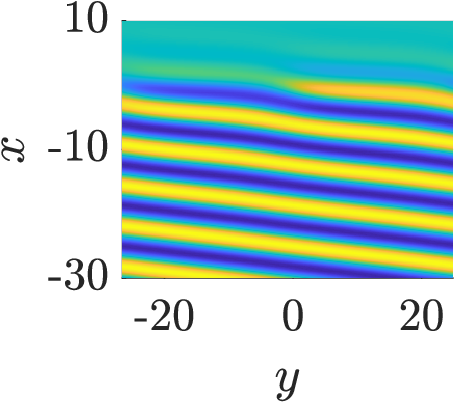}
         \includegraphics[trim=0 0 0 0,clip,width=0.8\textwidth]{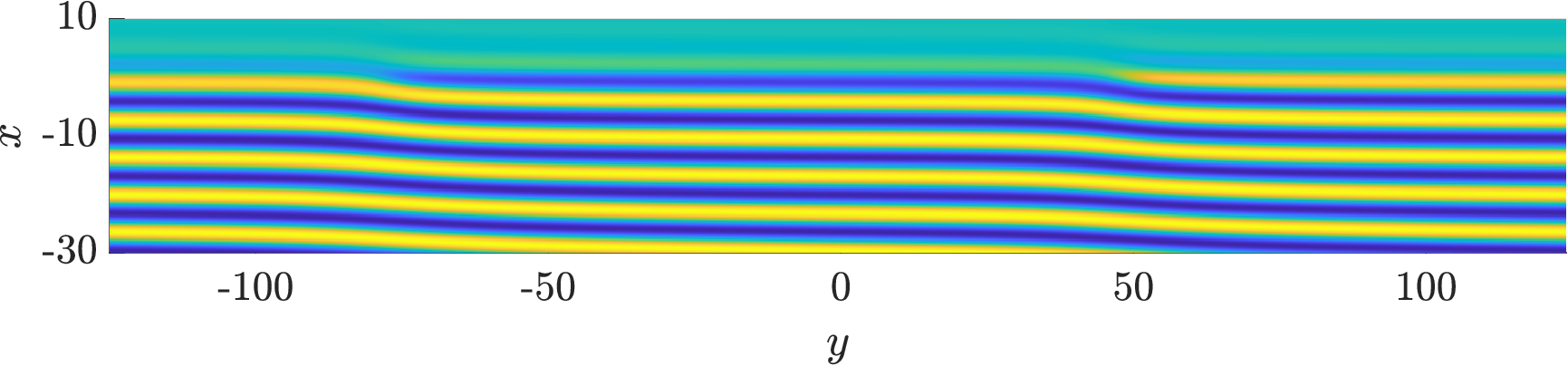}\\
    \includegraphics[trim=0 0 0 0,clip,width=0.49\textwidth]{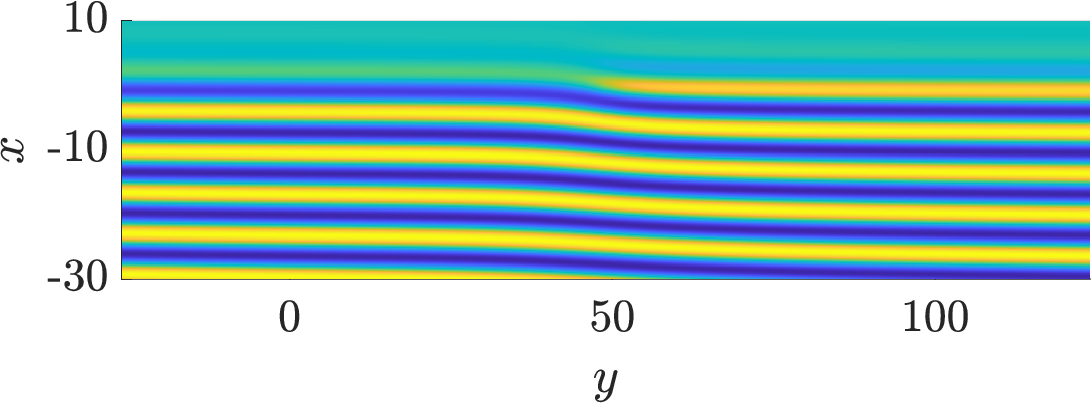}
    \includegraphics[trim=0 0 0 0,clip,width=0.49\textwidth]{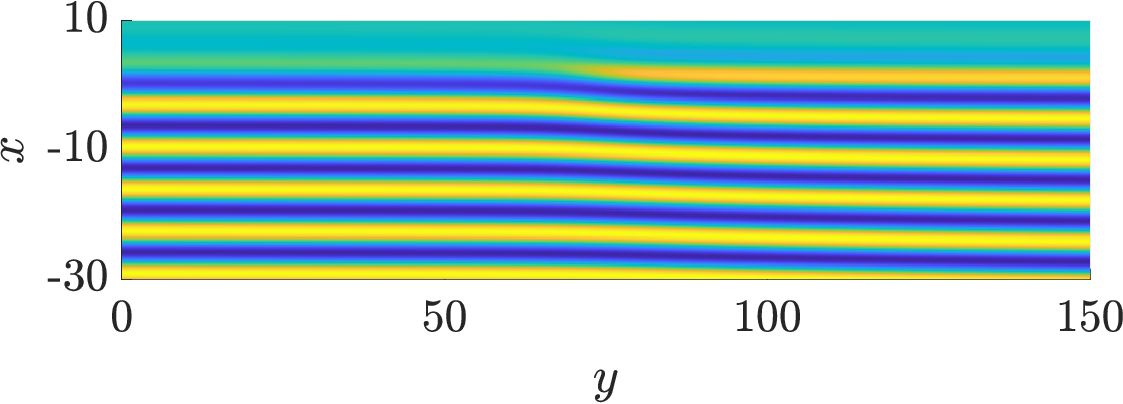}
     \caption{Plots of solutions of \eqref{e:aSH2}-\eqref{e:aSH4}  near quenching interface in original coordinates for $\tilde c_x = 10^{-3}$ fixed for a range of $\tilde k_y$ values:  $\tilde k_y = 118.23...$ (top left), $\tilde k_y = 25.13...$ (top right ). Bottom row illustrates delocalization of dislocation defect both in $x$ and $y$ for small $\tilde k_y$, with a zoom-in near a defect for $\tilde k_y =25.13...$ (left) $\tilde k_y = 4.35...$ (right). Note that the odd symmetry in Swift-Hohenberg creates two antisymmetric dislocation-type defects, a covering symmetry visible also in the phase-diffusion approximation through the dependence of the strain-displacement relation on $2\tilde{\phi}$, only.}\label{f:ash-sol}
\end{figure}

\paragraph{Comparisons between phase-diffusion and Swift-Hohenberg}
 Figure \ref{f:aSH} gives slices of the moduli space for \eqref{e:aSH} with $\tilde c_x$ fixed and shows that the surface is a graph $k_x = k_x(k_y,\tilde c_x)$ for $(k_y,\tilde c_x)\sim0$. Curves, which are plotted over the scaled wavenumber $\tilde k_y = k_y/\tilde c_x$, show good agreement with the heteroclinic asymptotics of Section \ref{s:origin}, with a transition around $k_y/\tilde c_x \sim 6$ between a localized defect near the quenching interface to the delocalized heteroclinic selecting smaller wavenumbers; see Figure \ref{f:ash-sol} for plots of relevant solutions.

%  switch

% \begin{figure}[h]
%     \includegraphics[trim=0 0 0 0,clip,width=0.19\textwidth]{figures/cx1e-3left.eps}
%          \includegraphics[trim=0 0 0 0,clip,width=0.8\textwidth]{figures/cx1e-3right.eps}\\
%     \includegraphics[trim=0 0 0 0,clip,width=0.99\textwidth]{figures/cx1e-3bottom.eps}
% 	    \includegraphics[trim=0 0 0 0,clip,width=0.99\textwidth]{figures/cx1e-3bottom_new.eps}
%      \caption{Plots of solutions of \eqref{e:aSH2}-\eqref{e:aSH4}  near quenching interface in original coordinates for $\tilde c_x = 10^{-3}$ fixed for a range of $\tilde k_y$ values:  $\tilde k_y = 118.23...$ (top left), $\tilde k_y = 25.13...$ (top right ), $\tilde k_y = 8.61...$ (middle row), $\tilde k_y = 5.39...$ (bottom row, only part of $x$-domain shown). Dislocation defect delocalizes both in $x$ and $y$ for small $\tilde k_y$ as can be seen in the bottom row graphics.
% }\label{f:ash-sol}
% \end{figure}

Varying the anisotropy coefficient $\beta$ and the parameter $\mu_0$, we also show how this phase transition depends on system parameters. As expected, the strength of non-adiabatic effects increases with $\mu_0$ as averaging is less effective, and the strain $1-k$ on the stripes created at small $k_y$ increases, roughly proportional to $\mu_0$ as predicted by the amplitude $\mu_0/16$ of the strain-displacement relation. The location of the transition appears to be roughly independent of $\mu_0$, in agreement with our derivation above. Varying the strength of anisotropy does affect the transition. Stronger anisotropy narrows the plateau where delocalized defects determine wavenumber selection. Very weak and in particular vanishing anisotropy lead to non-monotone dependence of $k$ on $k_y$ which is beyond the scope of this paper.

\section{Conclusions and discussion}\label{s:7}

We investigated directional growth of striped phases in the absence of instabilities and for weakly oblique orientation of stripes relative to the boundary. In a reduced phase-diffusion approximation, we established existence of simple, resonant growth mechanisms and derived universal asymptotics in limiting regimes. Our results compare well with computations in a Swift-Hohenberg equation where instabilities are suppressed by weak anisotropy. 

Many of our results can be rephrased in coarse terms. For parallel stripes, we had earlier found that very small speeds cause maximal strain, given by the minimum of the strain-dispersion relation, which decreases up to a dynamically averaged (harmonic average) strain for large speeds. Zero speeds and growth at larger angles yield zero strain, with selected wavenumber given by the (energy-minimizing) average of the strain-displacement relation. At small angles, $k_y\sim 0$, the growth process is mediated by the emergence of a point defect at the boundary, which undergoes a delocalization bifurcation at a critical value, similar in character to the codimension-two bifurcation from a hyperbolic homoclinic orbit to a saddle-node homoclinic orbit. The growth process is described well by a glide motion of the defect along the boundary of the patterned region, adding one stripe once the defect has moved by one period along the boundary. In our asymptotics, we identify the glide motion in the absence of growth, $c_x=0$, when a non-equilibrium strain $k\neq \dashint g$ is imposed in the far field: the nonequilibrium strain $k$ drives the defect at a finite speed $c_y(k)$, $c_y'\neq 0$.
Then, for a growth process with given speed $c_x$ and angle $k_y$, the selected wavenumber $k$ adjusts such that 
%The selected wavenumber $k$ in a growth process with given speed $c_x$ and angle $k_y$ then adjusts such that 
the induced glide speed $c_y(k)$ corresponds to compatible defect motion by one $y$-period $2\pi/k_y$ while one stripe is grown across the interface, in time $2\pi/(c_xk_x)$. The effective wavenumber used in the scaling, $k_y/c_x=k_x/c_y\sim 1/c_y$, is at leading order simply the inverse glide speed.  From this perspective, the $\tilde{k}_y$-dependent contribution to the strain stems from drag in the glide motion of the defect. The $c_x$-dependence can be understood as in \cite{beekie} as an interaction between dislocation over the finite distance $2\pi/k_y$,  leading to an effective deceleration of the glide motion and reduced strain.

\paragraph{Effect on energy densities} Our results can also be interpreted from an energetic point of view. The Swift-Hohenberg equation in the unquenched form is a gradient flow to the energy $\int E$ with energy density $E=\frac{1}{2}(\Delta u+u)^2 -\frac{\mu}{2}u^2+ \frac{1}{4}u^4$. Among the striped patterns there is a unique wavenumber $k_\mathrm{zz}=1+\rmO(\mu^2)$ that minimizes the energy per unit volume. The wavenumber $k_\mathrm{zz}$ happens to coincide with the onset of the zigzag instability in $k<k_\mathrm{zz}$ in the isotropic case, although this instability is suppressed in the anisotropic setting. Periodic patterns do in fact minimize the energy density in one space-dimension \cite{peletier} and one typically sees convergence to periodic patterns and energy densities vary close to the minimizer in large bounded one-dimensional domains. In higher dimensions, proofs that periodic patterns minimize energies do not appear to be available, and generic initial conditions do not converge to periodic patterns. Defects and boundary conditions play an important role both in the organization of stable stationary states and in the selection of wavenumbers. 

The present results demonstrate the effect of growth on energy in the bulk. The energy minimizing wavenumber corresponds to $k_x=1$ in the phase-diffusion approximation, such that the square deviation $(k_x-1)^2$ is a good approximation for the energy density of the pattern in the bulk, away from the interface. The selection of the energy minimizer at $k_y\neq 0,\ c_x=0$ echoes the selection of periodic patterns with minimal energy by grain boundaries \cite{lloydscheel}. At $k_y=c_x=0$, the boundary does not select a specific wavenumber but rather (significantly) narrows the band of compatible wavenumbers from $\rmO(\sqrt{\mu})$ to $\rmO(\mu)$, a mechanism also observed in point defects; see for instance  \cite[\S 4.4]{ssradial} for the case of a focus defect. For wavenumbers outside of the compatible band, one usually sees diffusive repair between the selected wavenumber and the imposed farfield wavenumber, as in the case of grain boundaries, or drift of phase and defects, as in the case $k_y\neq 0$, $c_x=0$.  

For nonzero speeds, the growth process selects a unique wavenumber away from the energy minimizer: The fact that $k_x<1$ guarantees that energy, inserted into the system at the moving quenching line,  is stored in the bulk at a constant density. In other words, the gradient dynamics are driven by a localized energy source and relax to equilibrium in the bulk away from the source, albeit not the energy-minimizing ``thermodynamic equilibrium''. Our results show that such a relaxation to the energy minimizer occurs only in the limit  $k_y\to\infty$, that is, for large angles between rolls and quenching line, or for vanishing non-adiabatic effects, $\kappa=0$ or $g\equiv const$. It does not seem obvious how one might quantify the stored energy in the system directly from energetic considerations at the quenching line.

For larger speeds, beyond the validity of the phase-diffusion approximation, one finds selected wavenumbers close to the wavenumbers selected by free invasion fronts \cite{GS1}. In a Ginzburg-Landau approximation, these select the minimium energy solutions. Higher-order corrections in the Swift-Hohenberg equation show however that the selected wavenumber does not correspond to the energy minimizer. In fact, most patterns created through directional quenching have wavenumbers below $k_\mathrm{zz}$ and are thus zigzag unstable in the isotropic case, although the instability may spread more slowly than patterns are created at the quenching line \cite{zigzag}. 

\paragraph{Other models of growth: heterogeneities and dynamic boundary conditions}
We also remark that several other growth processes also induce wavenumber selection phenomena which collapse the ``Busse Balloon" of possible wavenumbers supported in a homogeneous spatial domain \cite{Busse_1978}. For example, if the sharp quenching step with $c_x=0$ is replaced by a slowly varying parameter ramp, the band of compatible wavenumbers is significantly narrowed \cite{pomeauzaleski}. One could also model growth by restricting to a bounded, or semi-bounded domain with dynamic boundary. Various types of boundary conditions and their wavenumber selection properties in the wake were studied in the stationary case \cite{morrissey};  see also \cite{vanGorder21} and references therein for a review of other work in this direction. Motivated by precipitation and deposition phenomena, traveling source terms could also be used to force a system out of equilibrium and select wavenumbers in the wake \cite{thomas, dipstripe,gs5}.

\paragraph{Further directions: phenomena, theory, and experiments}
Looking forward, we hope that this glimpse into the role of point defects in growth of crystalline phases can be extended, including for instance the effect of zigzag instabilities associated with wrinkling. More mathematically, of the many phenomena described here, it would be interesting to analyze the heteroclinic bifurcation at the origin, finding in particular better asymptotics near the critical value of $\tilde{k}_y$. One may also hope to better understand some of the asymptotic expansions derived here, adding mathematical rigor, or relating them more directly to our understanding of dislocations, their farfield, and interaction properties.

 We hope that some of the predictions here can be confirmed in experiments; see \cite{electronematics} for a current overview of experimental setups in the context of  electroconvection with nematic liquid crystals. Approximation of dynamics by a Ginzburg-Landau equation has been confirmed quantitatively in many experiments, potentially allowing for quantitative comparisons with our results; see for instance \cite{fl} and references therein. A setup where applied currents can be controlled locally would then allow experiments that test some of our predictions. Most notably, it would be interesting to observe the non-monotonicity of strains in speed for small angles and compare the related dynamics of dislocation-type point defects near the quenching line and with the glide motion of free dislocations in the Ginzburg-Landau equation \cite{pismenmobility}.

% comment on experimental consequences ; lots of experiments described in ; quantitative comparisons with CGL described in \cite{fl}

\bibliographystyle{abbrv}
\bibliography{woq_rev}

\end{document}